\documentclass[12pt]{article}
  \usepackage{amsfonts,amsmath}
  \usepackage[square,sort&compress,comma,numbers]{natbib}
   \def\R{\mathbb{R}}
   \def\N{\mathbb{N}}
   
   \def\1{{\rm I\mskip -10.5mu 1}}
   
   \def\e{{\varepsilon}}
   \def\D{{\nabla}}
   \def\lo{\mathop{\longrightarrow}}

   \def\d{{\delta}}
   \def\l{{\lambda}}

   \def\cB{{\cal B}}
   \def\cC{{\cal C}}
   \def\cD{{\cal D}}
   
   \def\cF{{\cal F}}

   \def\cI{{\cal I}}
   \def\cK{{\cal K}}
   \def\cL{{\cal L}}
   
   \def\cN{{\cal N}}

   \def\cR{{\cal R}}
   \def\cS{{\cal S}}

   \def\supp{\mathop{\rm supp}\nolimits}
   
   \def\dist{\mathop{\rm dist}\nolimits}

   \def\loc{\mathop{\rm loc}\nolimits}

   \def\no{\noindent}
   \def\proof{\mbox {{\underline {\sf Proof}} \hspace{2mm}}}
   \def\qed{{\hfill {\em q.e.d.}\\\vspace{1mm}}}

   \newcommand{\beq}{\begin{equation}}
   \newcommand{\eeq}{\end{equation}}

\newtheorem{df}{Definition}[section]
\newtheorem{prop}[df]{Proposition}
\newtheorem{lemma}[df]{Lemma}
\newtheorem{teo}[df]{Theorem}
\newtheorem{rem}[df]{Remark}

\newtheorem{cor}[df]{Corollary}

 \newcommand{\sezione}[1]{\section{#1}\setcounter{equation}{0}}

\begin{document}

   \title{Infinitely many positive standing waves for Schr\"odinger
     equations with competing coefficients}
  \author{Giovanna Cerami\thanks{ Dipartimento di Matematica,
        Politecnico di Bari,  Via Amendola 126/B -  70126
        Bari, Italia.
        E-mail: {\sf  cerami@poliba.it}}
    \and
        Riccardo Molle\thanks{Dipartimento di
        Matematica, Universit\`a di Roma ``Tor Vergata'', Via
                della Ricerca
        Scientifica n$^o$ 1 - 00133 Roma. E-mail: {\sf
                  molle@mat.uniroma2.it}}
        }
\date{}

 \maketitle

{\small {\sc \noindent \ \ Abstract.} -
The paper deals with the equation $-\Delta u+a(x) u +b(x)u^q -u^p =
0$, $u \in H^1(\R^N)$, whith  $N\ge 2$, $1<q<p,\ p<{N+2\over N-2}$ if
$N\ge 3$, $\inf a>0$, $a(x)\to a_\infty$ and $b(x)\to 0$ as
$|x|\to\infty$.

When $a(x)\le a_\infty$ and $b(x) = 0$ only a finite number of
positive solutions to the problem is reasonably expected. 

Here we prove that the presence of a nonzero term $b(x)u^q $ with
$b(x)\geq 0, \ b(x)\neq 0,$ under suitable assumptions on the  decay
rates of $a$ and $b,$ allows to obtain  infinitely many 
positive solutions.

}

\vspace{3mm}

\sezione{Introduction}

In this paper we deal with the question of finding multiple positive
solutions to: 

$$(P) \hspace{1cm}\qquad\left\{\begin{array}{l}
-\Delta u+a(x) u +b(x)u^q -u^p = 0 \ \ \ \mbox{in} \ \ \R^N\\
u \in H^1(\R^N)
\end{array}\right.$$
where $N\ge 2$, $1<q<p$ and $p<2^*-1 = \frac{N+2}{N-2}$ if $N\ge 3$.
   The potential
$a(x)$ and the coefficient $b(x)$ are non negative functions, not
required to enjoy symmetry, such that 
$$(h_1) \hspace{1,5cm} \lim_{|x|\to +\infty} a(x) = a_\infty >0
\hspace{2cm} \lim_{|x|\to +\infty} b(x) = 0$$ 

Problems like $(P)$  have been widely investigated during last four
decades: it is well known that the interest  in studying Euclidean
scalar field equations has been motivated not only by their strong
relevance in Physics and Mathematical Physics  but also by their
mathematical features that make them challenging to the
researchers. Indeed, $(P)$ has a variational structure, but a lack of
compactness, due to the invariance of $\R^N$ under the action of the
non compact group of translations, prevents the use of  variational
methods in a standard way. 

Starting from the classical by now papers (\cite{BL,S}),
several existence and multiplicity results have been stated and 
qualitative properties  have been studied of the solutions to
equations $ -\Delta u+a(x) u  = f(x,u(x))$.  The earliest results
were obtained in radially symmetric situations, taking advantage of
the compact embedding in $L^p (\R^N),\ p\in (2, 2N/(N-2))$ of the
subspace of $H^1(\R^N)$ consisting of radial functions.  
Subsequently  many different devices have been exploited to face the
difficulties of  non symmetric situations and remarkable progresses
have been obtained also in this direction, even if many  questions
still remain open. 

Describing all the various and interesting contributions, without
forgetting something, is not an easy matter, so we refer interested
readers to a survey paper \cite{Cmilan} for
an introduction to the motivations for studying scalar field equations
and a description of the development of the research. 

Here, in order to  describe the purpose of our research, we need to
recall some results on $(P) $ when $b(x) = 0.$ Dealing with this
problem, one at once realizes that the topological framework is quite
different according the way in which $a(x)$ approaches its limit at
infinity. 
When $a(x)$ goes to $a_\infty$ from below the existence of a positive
ground state solution can be proven by a minimization method together
concentration compactness arguments (see \cite{PLL}). Conversely when $a(x)
> a_\infty$ it is not difficult to see that $(P)$ has not ground state
solutions and positive solutions, when they exist, do not correspond
to the infimum of the energy functional on the natural constraint, but
must be searched  by subtle topological and variational tools at
higher energy levels \cite{BLi,BPLL}.   

The same striking difference appears when one looks for multiple
solutions. Under  suitable slow decay assumptions on $a(x)$, in the
first case the existence of infinitely many \textit{changing sign
  solutions} has been shown in \cite{CdVS}, while in
the second case  the existence of infinitely many \textit{positive
  multi-bump solutions} has been proved first in \cite{CPS1} under an
additional smallness condition on the oscillation of $a(x)- a_\infty,$
subsequently  without that condition in \cite{CMP} for potentials
having ``dips'' then, in general,  in 
\cite{DvS2,DPWY} in the planar case and, finally, in any dimension of
the space $\R^N$ in \cite{MPlincei,MPpreprint}. 

It is worth remarking that, when $a(x) > a_\infty,$ the positive
multi-bump solutions rougly speaking  appear, when one tries to
maximize the energy functional between functions having a fixed number
of bumps,  as  the result of a delicate equilibrium between the
attractive effect of $a(x)$ and  the repulsive disposition of the
bumps (which can be tought as nothing but  solutions of the limit
problem $-\Delta u+a_\infty u  -u^p = 0$).   

The same heuristic argument can be used to understand why, when $a(x)<
a_\infty,$  one cannot expect to find infinitely many positive multi
bump solutions. Indeed, since the interaction between masses of the
same type makes the energy decrease, if one tries to maximize the
energy functional between multi-bump positive functions, the repulsive
effect of $a(x)$ plus the repulsive effect of the masses each other
push the masses at infinity,  on the contrary, if one attempts to
minimize the functional,  reasonably can hope to find just a finite
number of multi-bump solutions and not infinitely many,  because, when
the number of bumps increases, the  attraction due to  $a(x)$  added
to  the attractive disposition of the masses each other makes
impossible an equilibrium state of the bumps  and on the contrary
forces them to collapse.  

The object of our investigation is to understand how the presence of
another nonlinear term  can modify when $a(x)\leq a_\infty $ the above
described situation opening the possibility of getting infinitely many
positive solutions.  The result we obtain, which is contained in the
following theorem, states that this actually can happen: a ``small''
nonlinear perturbation having a coefficient whose  effect is
``competing'' with the attractive one of $a$ allows to find for any
$k\in \N$ a  positive solution having exactly $k$ bumps:

\begin{teo}
	\label{T}
	Let $a(x)$ and $b(x)$ satisfy $(h_1)$ and 
	
	$(h_2) \hspace{1,5cm} a_0:=\inf_{\R^N}a(x)>0 ; \ \  \alpha(x)
        =  a_\infty - a(x)\geq 0;$ 
	
	\vspace{0,2cm} 
	$(h_3) \hspace{1,5cm} b(x)\geq 0 ,\ \ \ b\in L^\infty(\R^N)$.
	
	\vspace{0,2cm} 
	Assume that $\eta\in(0,\sqrt{a_0})$ and $c>0$ exist such that 
	
	\vspace{0,2cm}
	$(h_4)	\hspace{1,5cm} \alpha(x)\le c\, e^{-\eta\, |x|},\ \ \
        \ \forall x\in\cD, \hspace{2cm} 
	\lim_{|x|\to+\infty,\ x\in \cD}b(x)\, e^{\eta\, |x|}=+\infty$
	
	\vspace{0,2cm} 
	\noindent where $\cD=\{r\theta\in\R^N$ : $r>0,\
        \theta\in\Theta\},$ $\Theta\subset S^{N-1}$ open, and that
        $B_d(\bar \zeta)\subset \overline \cD$ for some $\bar
        \zeta\in\Theta$ and $d\ge {1\over 2}.$ 
	
	Then $B>0$ exists such that if $|b|_{L^\infty}\,<\,B$  problem
        $(P)$ has infinitely many positive solutions. 
\end{teo}

Although  our main interest is devoted to the case in which $\alpha
(x)\not\equiv 0 $  we remark  that when  $\alpha( x) \equiv 0 $
Theorem \ref{T} gives a result comparable to that stated in \cite{CPS1}, but
under a weaker decay condition because  imposed in a conical region
instead of all the space.

Furthermore we stress the fact that some assumptions could be even
more weakened: $(h_4)$  could be substituted by the more technical
one:  
$$(h'_4) \ \ \  \exists\, \zeta_n\ \ s.t\ \ \  |\zeta_n | > 2n\ \ \
\mbox{and}\ \   \lim_{n\to \infty}\frac{\inf_{B_{n/2} (\zeta_n)} b}
{\sup_{B_{n/2} (\zeta_n)} \alpha}=\lim_{n\to \infty}\left(
  \inf_{B_{n/2} (\zeta_n)}b\right) e^{2n\eta} = +\infty$$ 
which asks a control of the asymptotic behaviour of $a$ and $b$ only
on a sequence of balls whose centres and radii go to
infinity. Moreover we observe that  the competing effect of $b$ is
mainly due to its behaviour at infinity, so $b$ could be allowed to be
negative outside the balls. We have chosen to state and prove Theorem
\ref{T} under assumptions $(h_3)$ and $(h_4)$ to avoid  more
technicalities. 

It is also worth observing that in our assumptions problem $(P)$ may
or may not have a positive ground state solution. In the appendix
examples are exhibited to show that both  cases occur.    

The method we use to prove theorem \ref{T} is purely variational and
inspired by that of \cite{CPS1}. 

The paper is organized as follows: in  section two some useful results
are recalled and the variational framework in which the problem is
studied is introduced, in section three a max-min procedure is
developed to find in suitable classes of multi bump positive functions
good candidates to be critical points, sections 4 and 5  contains the
proof that actually the functions found in section 3 are solutions
when $|b|_{L^\infty}$ is suitably small,  some examples are contained
in the appendix.

\sezione{Variational framework, useful facts and remarks, tools }

 Throughout the paper we make use of the following notation:
 
 \begin{itemize}
\item $H^{1}({\mathbb R}^{N})$ is the usual Sobolev space endowed with
  the standard scalar product and norm 
\begin{displaymath}
(u, v)=\int_{\mathbb R^N}[\nabla u \nabla v+a_\infty uv]dx;\qquad
||u||^{2}=\int_{{\mathbb R}^N}\left[|\nabla u|^{2}+a_\infty u^{2}\right]dx. 
\end{displaymath}
\item $H^{-1}$ denotes the dual space of $H^1(\mathbb R^N)$.
\item $L^q(\Omega)$, $1\leq q \leq +\infty$, $\Omega \subseteq \mathbb
  R^N$, denotes a Lebesgue space, the norm in $L^q(\Omega)$ is denoted
  by $|u|_{q, \Omega}$ when $\Omega$ is a proper subset of $\mathbb
  R^N$, by $|\cdot|_p$ when $\Omega=\mathbb R^N$. 
\item For any $\rho>0$ and for any $z\in \mathbb R^N$, $B_\rho(z)$
  denotes the ball of radius $\rho$ centered at $z,$ \and for any
  measurable set $ \mathcal{O} \subset \R^N, \ |\mathcal{O}|$ denotes
  its Lebesgue measure. 

\item $c,c', C, C', C_i$ denote various positive constants.
      \\
\end{itemize}

 Problem $(P)$ is variational, its solutions can be searched as critical
 point of the ``action''  functional $I:H^1(\R^N)\to \R$
$$
I(u)={1\over 2}\int_{\R^N}(|\D u|^2+a(x)u^2)dx
+{1\over q+1}\int_{\R^N}b(x)|u|^{q+1}
-{1\over p+1}\int_{\R^N}|u|^{p+1}.
$$

In what follows we need also to consider the so called ``limit''
problem of $(P)$: 
$$
(P_\infty)
\qquad\left\{\begin{array}{l}
-\Delta u+a_\infty u=|u|^{p-1}u\\
u\in H^1(\R^N)
\end{array}\right .
$$
and the related functional
$I_\infty:H^1(\R^N)\to \R$ defined by
$$
I_\infty(u)={1\over 2}\int_{\R^N}(|\D u|^2+a_\infty u^2)dx-{1\over
  p+1}\int_{\R^N}|u|^{p+1}dx.
$$

The following lemmas are well known and contain some useful information on
the solutions of  $(P_\infty)$  (see f.i.
\cite{Cmilan} and \cite{AM}, respectively).

\begin{lemma}
\label{N2}
 Problem $(P_\infty)$
has a positive ground state solution $w$, unique up to translation,
 radially symmetric,
decreasing when the radial coordinate increases and such that
\begin{equation}
\label{N2*}
\lim_{|x|\to + \infty}|D^j w (x)||x|^{\frac{N-1}{2}}e^{\sqrt{a_\infty}|x|}=d_j>0,
\hspace{3mm} d_j\in\R,\ j=0,1.
\end{equation}
\end{lemma}

\begin{lemma}
	\label{N2'}
	Set, for all $y \in \R^N,$
	\beq
	\label{N2'.1}
	w_y (x) = w(x-y),
	\eeq
	and
	$$
	\Phi = \left\{ w_y : \ y\in\R^N\right\}.
	$$
	
	Then $\Phi$ is nondegenerate, namely the following properties are true:
	
	i)\ $(I_{\infty})''(w_y)$ is an index zero Fredholm map, for
        all $y \in \R^N$; 
	
	ii) Ker $(I_{\infty})''(w_y)$ = span $\left\{ \frac{\partial
            w_y}{\partial x_j}: 1 \leq j \leq N \right\} =
        T_{w_y}\Phi,$ 
	
	$T_{w_y}\Phi$ being the tangent space to $\Phi$ at $w_y.$
\end{lemma}

We remark also that, setting
$$m_{\infty} := I_{\infty} (w) ,$$
$m_{\infty}$ can be characterized as:
\beq
\label{2.1}
m_{\infty} = \min_{\gamma \in \Gamma}\max_{t\in [0,1]} I_{\infty}(\gamma (t)),
\eeq
where, denoting by $e\in H^1( \R^N)$ any point for which $I_{\infty}(e)<0$,
\beq
\label{Gamma}
 \Gamma =\{\gamma \in C([0,1], H^1( \R^N))\  :\  \gamma(0)=0,\
   \gamma(1)=e\}.
\eeq

Let $\delta\in(0,1)$ be chosen so small that
\beq
\label{delta}
\begin{array}{ll}
\vspace{2mm}
\mbox{\em (i)}\ a_0-p\d^{p-1}>0, &\quad \mbox{\em (ii)}\ {a_0 \over
    2}-{\d^{p-1}\over
    p+1}>0,\\ 
\mbox{\em (iii)}\ a_0-\d^{p-1}>\eta^2, 
& \quad \mbox{\em (iv)}\  a_0\d-2^{p-1}\d^p>0.
\end{array}
\eeq

According to the choice of $\delta$ we fix $R$ so large that
$w(x)<\delta$ for $|x|> R/2$.

For all function $u\in H^1(\R^N)$, $u\ge 0$, we
use the notation
$$
u_\d=\min\{u,\d\},\qquad \mbox{ and }\quad u^\d=\max\{0,u-\d\}
$$
and we call $u_\d$ and $u^\d$ respectively the submerged part under
$\delta$ and the emerging part above $\d$ of $u$. Clearly
$$
u=u_\d+u^\d.
$$
Fixing $\d$ and $R$ as before indicated and $k\in\N\setminus\{0\}$ we
set
$$
\cK_k=\ \left\{\begin{array}{lc}
\R^N & \mbox{ if }k=1\\
\{(x_1,\ldots,x_k)\in(\R^N)^k\ :\ |x_i-x_j|\ge 2R,\ i\neq
j\} & \mbox{ if }k>1
\end{array} \right.
$$
and if $(x_1,\ldots,x_k)\in\cK_k$ we say that a function $u\in
H^1(\R^N)$ is emerging around the points $x_1,\ldots,x_k$ in balls of
radius $R$ if $u\ge 0$ and 
$$
u^\delta=\sum_{i=1}^k u_i^\delta,\mbox{
where
}
u_i^\delta\in H^1_0(B_{R}(x_i)),\ u_i^\delta\ge 0,\
u_i^\delta\not\equiv 0,\ i\in\{1,\ldots,k\}.
$$

We are now able to introduce the classes of functions in which we
look for solutions, indeed for all $(x_1,\ldots, x_k)\in\cK_k$ we set
\begin{eqnarray*}
S_{x_1,\ldots,x_k}\hspace{-2mm} &=\{&\hspace{-2mm}
u\in H^1(\R^N)\  : u \geq  0, \  u \mbox{
   emerging around }(x_1,\ldots,x_k),\\
& &
 I'(u)[u_i^\delta]=0,\   \beta_{x_i}(u)=0 \ \forall i\in\{1,\ldots,k\}\},
\end{eqnarray*}
where 
$$
\beta_{x_i}(u)\ =\ \frac{1}{\int_{B_{R}(x_i)}(u_i^\delta)^2}
\int_{B_{R }(x_i)} 
 (u_i^\d(x))^2(x-x_i)\, dx
$$
is a kind of ``local'' barycenter.

It is useful for what follows to introduce the functionals $J$ and
$J_\infty$ defined on the set $\{v\in H^1(\R^N)\ :\ |\supp
v|<+\infty\}$ respectively by
 \begin{eqnarray*}
J(v) & := & {1\over 2}\int_{\cS_v} (|\D v|^2+a(x)v^2)dx+\delta\int_{\cS_v} a(x) v\,dx\\
 & & +{1\over q+1}\int_{\cS_v} b(x)(\delta+v)^{q+1}dx-{1\over
   p+1}\int_{\cS_v}(\delta+v)^{p+1}dx\\
& & -{\delta^{q+1}\over q+1}\int_{\cS_v} b(x)\, dx +{\d^{p+1}\over p+1}\, |{\cS_v}|.
\end{eqnarray*}
\begin{eqnarray*}
J_\infty(v) & := & {1\over 2}\int_{\cS_v}(|\D v|^2+a_\infty v^2)dx+\delta\int_{\cS_v} a_\infty v\,dx\\
 & & -{1\over
   p+1}\int_{\cS_v}(\delta+v)^{p+1}dx +{\d^{p+1}\over p+1}\, |{\cS_v}|
\end{eqnarray*}
where ${\cS_v}:=\supp v$.

We remark that for all $u\in H^1(\R^N)$
\beq
\label{*_1}
I(u)=I(u_\delta)+J(u^\delta)
\eeq
\beq
\label{*_2}
I_\infty(u)=I_\infty(u_\delta)+J_\infty(u^\delta).
\eeq

Next lemmas describe the nature of the local non-smooth constraints
$I'(u)[u_i^\d]=0$ imposed to the functions belonging to sets
$S_{x_1,\ldots,x_k}$; as a consequence we also deduce
$S_{x_1,\ldots,x_k}\neq\emptyset$, for all $(x_1,\ldots,x_k)\in \cK_k$
and for all $k\in\N\setminus\{0\}$.

\begin{lemma}
\label{lMax}
Let $u\in H^1(\R^N)$, $u\ge 0$ be such that $u^\d\neq 0$.
Let assumptions $(h_1)$, $(h_2)$, $(h_3)$  be satisfied.
Then a number $B_1\in\R^+\setminus\{0\}$ exists so that, if the
condition $|b|_\infty<B_1$ is fulfilled, the function
$g:[0,+\infty)\to\R$ defined as
$$
g(t)=I(u_\d+tu^\d)
$$
has a unique maximum point in $(0,+\infty)$.
\end{lemma}

\noindent \proof
Setting $\cS=\supp u^\d$ and using (\ref{*_1}) we get
\beq
\label{B1}
\begin{array}{rcl}
\vspace{2mm}
g'(t)={d\over dt}\, J(tu^\d) & = &\left[\int_\cS (|\D u^\delta|^2+a(x)(u^\delta)^2)dx\right]
t+\int_\cS a(x)\d u^\delta dx \\
& &
+\int_\cS b(x)(\d+t u^\delta)^qu^\delta dx
- \int_\cS (\d+t u^\delta)^pu^\delta dx.
\end{array}
\eeq
By (\ref{delta})({\em iii}) we then obtain 
\beq
\label{2.6}
\begin{array}{rcl}
\vspace{2mm}
g'(0)=J'(0)&=&\int_\cS(a(x)+b(x)\delta^{q-1}-\delta^{p-1})\delta\,
u^\delta dx
\\
&\ge&
\int_\cS(a_0-\d^{p-1})\d \, u^\delta dx >0,
\end{array}
\eeq
and, since $\lim\limits_{t\to +\infty} g'(t)=-\infty$, we deduce that
$g$ must have a zero in $(0,+\infty)$.

Moreover
$$
g'''(t)=q(q-1)\int_\cS b(x)(\d+t u^\delta)^{q-2}(u^\delta)^3 dx
-p(p-1) \int_\cS (\d+t u^\delta)^{p-2}(u^\delta)^3dx.
$$
Now, since $q<p$, a constant $\bar c>0$, depending on $\d$, can be
found so that $c (\d+s)^{q-2}\le (\d+s)^{p-2}$, $\forall s\ge 0$ and
$\forall c\le\bar c$. 
Therefore, if $|b|_\infty$ is small enough $g'''(t)<0$, for all $t\in
[0,+\infty)$, so we conclude that $g'$ is concave and, hence, its zero
is unique. 

\qed

\begin{cor}
\label{C2.5}
Let assumptions $(h_1),\, (h_2),\, (h_3)$   be satisfied. 
Let $k\in\N\setminus\{0\}$, $(x_1,\ldots,$ $ x_k)\in\cK_k$ and $u\in
H^1(\R^N)$ be a function emerging around $(x_1,\ldots,x_k)$.
Then the same conclusion of Lemma \ref{lMax} holds true for the
functions
$$
g_j(t)=I\left(u_\delta+\sum_{i\neq j}u_i^\delta+tu_j^\delta\right)\qquad
j\in\{1,\ldots,k\}.
$$
\end{cor}

\noindent\proof
The argument is quite analogous to that of Lemma \ref{lMax} once we
write for all $j\in\{1,\ldots,k\}$
$$
g_j(t)=
I\left(u_\delta+\sum_{i\neq j}u_i^\delta\right)+J(tu_j^\delta).
$$
\qed

In what follows we denote by $\theta(u)\in (0,+\infty)$ the maximum point
of the function $g(t)$ and by $\theta_j(u)\in(0,+\infty)$ the maximum
point of $g_j(t)$, $j\in\{1,\ldots,k\}$.

\begin{rem}
\label{5*}
{\em It is now clear that, if $u$ is emerging around
  $(x_1,\ldots,x_k)$ and satisfies for all $j\in\{1,\ldots,k\}$ the
  conditions $I'(u)[u^\d_j]=0$, then $\theta_j(u)=1$ for all
  $j\in\{1,\ldots,k\}$, furthermore $I(u)=\max\limits_{t\ge 0}I(u_\d+\sum_{i\neq j}u_i^\delta
+tu_j^\delta)$ and $J(u^\d_j)=\max\limits_{t\ge 0} J(tu^\d_j)$ for all
$j\in\{1,\ldots,k\}$.
}
\end{rem}
 
\begin{df}
\label{5**}
Let $u\in H^1(\R^N)$ be a function emerging around $(x_1,\ldots,
x_k)\in\cK_k$, we call $u_\d+\theta(u)u^\d$ the {\em projection of $u$
  on the set $\{u\in H^1(\R^N) \ :\ I'(u)[u^\d]=0\}$,
  $u_\d+\sum_{i\neq j} u_i^\d+\theta_j(u)u^\d_j$ the projection of $u$
  on $\{u\in H^1(\R^N)\ :\ I'(u)[u^\d_j]=0\}$. }
Moreover, if $u$ is such that $\beta_{x_i}(u)=0$ for all
$i\in\{1,\ldots,k\}$, the function
$u_\d+\sum_{i=1}^k\theta_i(u)u^\d_i$ is said  {\em projection of $u$
  on $S_{x_1,\ldots,x_k}$}.
\end{df}

\begin{lemma}
\label{L2.6}
Let assumptions $(h_1),\, (h_2),\, (h_3)$  be satisfied and $|b|_\infty<B_1$.
Then, for all $k\in\N\setminus\{0\}$, for all
$(x_1,\ldots,x_k)\in\cK_k$, $S_{x_1,\ldots,x_k}\neq\emptyset$.
\end{lemma}

\no \proof Let $\phi\in\cC_c(B_{R}(0))$ be a radially symmetric (around
the origin) function such that $\max_{B_R(0)}\phi>\d$.
Setting $u=\sum_{i=1}^k\phi(\cdot -x_i)$, the function $\bar
u=u_\d+\sum_{i=1}^k\theta_i(u)u^\d_i$ belongs to $S_{x_1,\ldots,x_k}$
because by definition of barycenters $\beta_{x_i}(\bar
u)=\beta_{x_i}(u)=0$.

\qed

Assume that $\delta$ satisfies (\ref{delta}) and set
$$
\cC=\{u\in H^1(\R^N)\ :\ 0\le u\le\d\ \mbox{ a.e. in }\R^N\}.
$$
\begin{lemma}
\label{L2.7}
Let assumptions of Lemma \ref{L2.6} hold.
Then 
\begin{itemize}
\item[{ A)}] $I$ is coercive and convex, hence weakly lower
  semicontinuous on $\cC$;
\item[{ B)}] $\forall k\in\N\setminus\{0\}$, $\forall
  (x_1,\ldots,x_k)\in\cK_k$,
$$
u\in S_{x_1,\ldots,x_k}\ \Rightarrow\ \begin{array}{rcc}
i) & J(u_i^\d)>0 & \forall i\in\{1,\ldots,k\} \\
ii) & I(u)>0. &
\end{array}
$$
\end{itemize}
\end{lemma}

\no \proof
{\em A)} \quad Let $u\in\cC$, by using (\ref{delta}){\em (ii)} we get 
\begin{eqnarray*}
\vspace{2mm}
I(u) & \ge & {1\over 2}\int_{\R^N}(|\D u|^2+a_0u^2)dx-{1\over p+1}\int_{\R^N}u^{p+1}dx
\\ \vspace{2mm}
& \ge & {1\over 2}\int_{\R^N}|\D u|^2 dx+\left({a_0 \over
    2}-{\d^{p-1}\over
    p+1}\right)\int_{\R^N}u^2dx\\
&\ge & c \|u\|^2.
\end{eqnarray*}
To show that $I$ is convex let us fix $u_1,u_2\in\cC$ and set
$h(t)=I(u_1+t(u_2-u_1))$. 
Thanks to   (\ref{delta}){\em (i)} we then obtain
\begin{eqnarray*}
\vspace{2mm}
h''(t)&=&\int_{\R^N}(|\D (u_2-u_1)|^2+a(x)(u_2-u_1)^2)dx
\\ \vspace{2mm}
 & & +\int_{\R^N}[qb(x)(u_1+t(u_2-u_1))^{q-1}-p(u_1+t(u_2-u_1))^{p-1}](u_2-u_1)^2dx
\\ \vspace{2mm}
& \ge & \int_{\R^N}|\D (u_2-u_1)|^2dx+\int_{\R^N}(a_0- p \d ^{p-1})(u_2-u_1)^2dx
\\ \vspace{2mm}
& \ge & 0
\end{eqnarray*}

Relation $(B)(i)$ follows straightly from the proof of Lemma
\ref{lMax} and Remark \ref{5*}, indeed we have
$$
J(u^\d_i)=\max_{t\ge 0} J(tu^\d_i)>J(0)=0.
$$

{\em (B)(ii)} is a consequence of {\em (A), (B)(i)} and (\ref{*_1}).

\qed

\begin{rem}
\label{7*}
{\em
We point out that when $a(x)\equiv a_\infty$ and $b\equiv 0$, the definition of
$S_{x_1,\ldots,x_k}$ still makes sense and Lemmas \ref{lMax},
\ref{L2.6}, \ref{L2.7}, Corollary \ref{C2.5}, hold.
In this case we use the notation $S^\infty_{x_1,\ldots,x_k}$.
Moreover, we can rephrase in this case too Remark \ref{5*} and
Definition \ref{5**} denoting by $u_\d+\theta^\infty(u)u^\d$ and by
$u_\d+\sum_{i\neq j}u^\d_i+\theta^\infty_j(u)u^\d_j$ the projections
of $u$ respectively on the sets $\{u\in H^1(\R^N)\ :\
(I^\infty)'(u)[u^\d]=0\}$, $\{u\in H^1(\R^N)\ :\
(I^\infty)'(u)[u^\d_j]=0\}$.
}
\end{rem}

Next lemma gives useful information on the set of the projections
of the ground state solutions of $(P_\infty)$ family.

\begin{lemma}
\label{L7*}
Let assumptions of Lemma \ref{L2.6} be satisfied.
Then $\{\theta(w_y)\ :\ y\in\R^N\}$ is a bounded set.
Moreover, for any sequence $\{y_n\}_n$, $y_n\in\R^N$, such that
$\lim_{n\to + \infty}|y_n|=+\infty$ the relation
$\lim_{n\to + \infty}\theta(w_{y_n})=1$ holds.
\end{lemma}

\no\proof
If for some sequence $\{y_n\}_n$,
$\lim_{n\to + \infty}\theta(w_{y_n})=+\infty$ would be true, then the
relation 
\beq
\label{pr2}
\begin{array}{rcl}
\vspace{2mm}
\lim\limits_{n\to+\infty} \theta (w_{y_n}) \left[ \int_{B_{R}(y_n)} (|\nabla w_{y_n}^\d|^2 +
a (x)(w_{y_n}^\d)^2)dx \right] + \d
\int_{B_{R}(y_n)}a(x) (w_{y_n})^\d dx  & &
\\
 +\int_{B_{R}(y_n)}b(x)(\d + \theta(w_{y_n})w_{y_n}^\d)^q w_{y_n}^\d
dx -\int_{B_{R}(y_n)}(\d + \theta(w_{y_n})w_{y_n}^\d)^p
w_{y_n}^\d dx&= & - \infty
\end{array}
\eeq
should be verified, contradicting
\beq
\label{pr}
\begin{array}{rcl}
\vspace{2mm}
 \theta (w_{\bar x}) \left[ \int_{B_{R}(\bar x)} (|\nabla w_{\bar x}^\d|^2 +
a (x)(w_{\bar x}^\d)^2)dx \right] + \d
\int_{B_{R}(\bar x)}a(x) (w_{\bar x})^\d dx  & &
\\
 +\int_{B_{R}(\bar x)}b(x)(\d + \theta(w_{\bar x})w_{\bar x}^\d)^q w_{\bar x}^\d
dx -\int_{B_{R}(\bar x)}(\d + \theta(w_{\bar x})w_{\bar x}^\d)^p
w_{\bar x}^\d dx&= & 0
\end{array}
\eeq
which, by definition, holds for all $\bar x\in\R^N$.

Now, let $\{y_n\}_n$ be a sequence of points in $\R^N$ such that
$\lim_{n\to + \infty}|y_n|=+\infty$.
Considering that $a(y_n)\longrightarrow a_{\infty}$,
$b(y_n)\longrightarrow 0$, as $n\to + \infty$, and that, up to a
subsequence, $\lim_{n\to + \infty}\theta(w_{y_n})=\bar\theta\in\R$,
(\ref{pr}) gives, as $n\to + \infty$, 
\beq
\label{pr1} 
\bar \theta \left[ \int_{B_{R}(0)} (|\nabla w^\d|^2
+ a_\infty(w^\d)^2)dx \right] + \d \int_{B_{R}(0)}a_\infty w^\d
dx -\int_{B_{R}(0)}(\d + \bar \theta w^\d)^p w^\d dx=0
\eeq 
from which we deduce $\bar\theta =1$, by Lemma \ref{N2} and the
uniqueness of $\theta^\infty(w)$.

\qed

We close this section by a lemma on the asymptotic decay of solutions of some
 variational inequalities. Its proof can be found in 
 \cite[Lemma 3.3]{CPS1}.

\begin{lemma}
\label{N*} Let $\cD \subset \R^N$ be closed, let  $\lambda, s \in
\R^+\setminus \left\{0\right\}$ and assume that $u\in H^1(\R^N)$
verifies
$$
 \hspace{1cm} \left\{
        \begin{array}{ll}
         -\Delta u+ \lambda u \leq 0         & \mbox{ in }\ \R^N \setminus \cD,\\
          0<u \leq s                     & \mbox{ in }\ \R^N \setminus \cD.
        \end{array}
\right.
$$
Then, for all $d \in (0, \sqrt{\lambda}), $ there exists a positive
constant $ C_d$ (depending on $\lambda , d, N$) such that
$$
u(x)\leq C_d s\, e^{-d\,\dist(x,\cD)} \ \ \ \ \ \ \ \ \ \ \  \forall x \in \R^N \setminus \cD.
$$
\end{lemma}

\sezione{Max-min scheme}

The purpose of this section is to find for all $k\in\N\setminus\{0\}$
a function emerging around $k$ points having the appropriate features
to be a critical point of the functional $I$.

From now on we assume that $(h_1),\, (h_2),\, (h_3)$   hold, and that
$|b|_\infty<B_1$ where $B_1$ is the number whose existence is stated
in Lemma \ref{lMax} and Corollary \ref{C2.5}.

By {\em (B)} of Lemma \ref{L2.7} we know that for all
$(x_1,\ldots,x_k)\in\cK_k$, $\inf_{S_{x_1,\ldots,x_k}}I\ge 0$,
  therefore we can define
\beq
\label{mu}
 \mu(x_1,\ldots,x_k):=\inf_{S_{x_1,\ldots,x_k}}I.
\eeq

\begin{prop}
\label{Pmin}
For all $k\in\N\setminus\{0\}$, for all
$(x_1,\ldots,x_k)\in\cK_k$, there exists $\bar u\in
S_{x_1,\ldots,x_k}$ such that
$$
I(\bar u)=\mu(x_1,\ldots,x_k).
$$
\end{prop}

\no\proof
Let $k\in\N\setminus\{0\}$ and $(x_1,\ldots,x_k)\in\cK_k$ be fixed and
let $\{u_n\}_n$, $u_n\in S_{x_1,\ldots,x_k}$  be a sequence such that
\beq
\label{9*}
\lim_{n\to + \infty}I(u_n)=\mu(x_1,\ldots,x_k).
\eeq
Then, using {\em (A)} of Lemma \ref{L2.7}, we obtain 
$$
c_1\|(u_n)_\d\|^2\le I((u_n)_\d)\le I(u_n)\le c_2\qquad
c_1,c_2\in\R^+\setminus\{0\}.
$$
We claim that $\left\{\|(u_n^\d)_i\|/|(u^\d_n)_i|_{p+1}\right\}$ is
bounded too, for all $i\in\{1,\ldots,k\}$. 
Otherwise, for some $i$ up to a subsequence, $\lim\limits_{n\to +
  \infty}{\|(u^\d_n)_i\|\over |(u^\d_n)_i|_{p+1}}=+\infty$ then   
\beq 
\label{8.-1}
\begin{array}{rcl}
\vspace{2mm}
\lim\limits_{n\to  + \infty}J\left({(u_n^\d)_i\over
    |(u_n^\d)_i|_{p+1}}\right)&\ge&\lim\limits_{n\to
   + \infty}\left[c_1{\|(u_n^\d)_i\|^2
\over|(u_n^\d)_i|_{p+1}^2}-c_2-\int_{B_{R}(x_i)}\left(\d+{(u_n^\d)_i\over
    |(u_n^\d)_i|_{p+1}}\right)^{p+1}dx\right]
\\
&=&+\infty
\end{array}
\eeq
which implies, in view of Corollary \ref{C2.5}, Remark \ref{5*}, and
Lemma \ref{L2.7}
\beq
\label{8.0}
\begin{array}{rcl}
\vspace{2mm} \lim\limits_{n\to +\infty} I(u_n)&=&\lim\limits_{n\to
+\infty}\left( I\left((u_n)_\d+\sum_{j\neq
i}(u_n^\d)_i\right)+\max_{t\ge 0}J(t(u_n^\d)_i)\right)
\\ \vspace{2mm}
& \ge & \lim\limits_{n\to+\infty}\left(I\left((u_n)_\d+\sum_{j\neq
i}(u_n^\d)_i\right)+J\left((u_n^\d)_i\over
    |(u_n^\d)_i|_{p+1}\right)\right)
\\
& = & +\infty,
\end{array}
\eeq 
contradicting (\ref{9*}).

Now, setting 
$$
\check u_n=(u_n)_\d+\sum_{i=1}^k{(u_n^\d)_i\over
  |(u_n^\d)_i|_{p+1}}
$$
we can assert the existence of $\check u\in H^1(\R^N)$ so that, up to
a subsequence
$$
\check u_n\longrightarrow \check u
$$
weakly  in $H^1(\R^N)$ and in $L^{p+1}(\R^N)$, a.e. in $\R^N$, strongly in
$L^2_{\loc}(\R^N)$ and $L^{p+1}_{\loc}(\R^N)$.
Therefore, for all $i\in\{1,\ldots,k\}$, $\check u_n\to \check u$
strongly in $L^{p+1}(B_R(x_i))$ and in $L^{2}(B_R(x_i))$, $0\le \check
u\le \d$ in $\R^N\setminus\cup_{i=1}^kB_R(x_i)$ and $\beta_{x_1}(\check
u)=0$.
Furthermore, by {\em (A)} of Lemma \ref{L2.7}
\beq
\label{8.2}
I(\check u_\d)\le\liminf_{n\to  + \infty}I((u_n)_\d).
\eeq
Lastly, let us define $\bar u=\check u_\d+\sum_{i=1}^k \theta
_i(\check u)\check u^\d_i$; then $\bar u\in S_{x_1,\ldots,x_k}$ and
it is the minimizer we are looking for.
Indeed

\begin{eqnarray*}
\vspace{2mm}
0<I(\bar u) & = & I\left( \check u_\d+\sum_{i=1}^k\theta_i(\check
u)\check u_i^\d\right)
\\ \vspace{2mm}
&\le & \liminf_{n\to +\infty} I\left( (u_n)_\d+\sum_{i=1}^k\theta_i(\check
u)(u_n^\d)_i\right)
\\ \vspace{2mm}
 &\le & \lim_{n\to +\infty} I\left( (u_n)_\d+\sum_{i=1}^k(u_n^\d)_i\right)
\\
&=&\mu(x_1,\ldots,x_k).
\end{eqnarray*}

\qed

\begin{rem}
{\em
We stress the fact that Proposition  \ref{Pmin} holds true when
$\alpha(x)=b(x)=0$.
Therefore, we define for all $k\in\N\setminus\{0\}$ and for all
$(x_1,\ldots,x_k)\in \cK_k$, 
$$
\mu^\infty (x_1,\ldots,x_k):=\inf_{S^\infty_{x_1,\ldots,x_k}}
  I^\infty(u)
$$
and we can assert this infimum is a positive minimum.
}\end{rem}

In what follows   
for all $(x_1,\ldots,x_k)\in\cK_k$, we set
$$
M_{x_1,\ldots,x_k}=\{u\in S_{x_1,\ldots,x_k}\ :\
I(u)=\mu(x_1,\ldots,x_k)\},
$$ 
$$
M^\infty_{x_1,\ldots,x_k}=\{u\in S^\infty_{x_1,\ldots,x_k}\ :\
I^\infty(u)=\mu^\infty(x_1,\ldots,x_k)\}.
$$

Next propositions describe some important features of any function
that realizes $\mu(x_1,\ldots,x_k)$.
First one describes the asymptotic behaviour of the submerged part of
any such function

\begin{prop}
\label{P3.3} Let $(x_1,\ldots,x_k)\in\cK_k$ and $u\in
M_{x_1,\ldots,x_k}$ then $u_\d$ satisfies \beq \label{eqs}
\left\{\begin{array}{ll}
-\Delta u+a(x)u+b(x)u^q=u^p &\mbox{ in }\R^N\setminus\supp u^\d\\
u=\d &\mbox{ on }\supp u^\d\\
u>0 &\mbox{ in }\R^N.
\end{array}\right.
\eeq
Moreover, let us fix $\eta_s\in(\eta,\sqrt{a_0-\d^{p-1}})$, then there
exists a constant $c>0$ such that 
\beq
\label{as}
0<u_\d<c\, e^{-\eta_s\dist(x,\supp u^\d)}.
\eeq
\end{prop}

\no\proof
In view of {\em (A)} of Lemma \ref{L2.7}, $I$ is coercive and convex
on the convex set
$$
\cL:=\{v\in H^1(\R^N)\ :\ 0\le v\le\d,\ v=\d\ \mbox{ on }\supp u^\d\}
$$
and $u_\d\in\cL$ is the unique, positive minimizer for the minimization
problem $\min\{I(v)\ :\ v\in\cL\}$, that is, $u_\d$ satisfies
$I'(u_\d)[v-u_\d]\ge 0$ for all $v\in\cL$.
This implies $-\Delta u+a(x)u+b(x)u^q=u^p$ for all $x\in\R^N$ such
that $u_\d<\d$.
On the other hand, by the choice of $\d$, $\bar u(x)=\d$ is a strict
supersolution of (\ref{eqs}) so the relation $u_\d<\d$ holds true for
all $x\in\R^N\setminus\supp u^\d$.
Relation (\ref{as}) follows straightly from Lemma \ref{N*}, by the
choice of $\d$, and by observing that since $u_\d$ solves (\ref{eqs}),
it also solves
$$
\left\{\begin{array}{lc}
-\Delta u+(a_0-\d^{p-1})u\le 0 &\mbox{ in }\R^N\setminus\supp u^\d\\
0<u_\d\le\d &\mbox{ in }\R^N\setminus\supp u^\d.
\end{array}\right.
$$

\qed

The following proposition states the conditions the gradient of $I(u)$,
subject to the $k$ local constraints coming from the definition of $S_{x_1,\ldots,x_k}$, has to verify when $u$ is a
minimizer on that class of functions.
The proof of it can be obtained arguing exactly as in \cite[Proposition
3.5]{CMP}.

\begin{prop}
\label{CPS3.6}
Let $k\in\N$, $(x_1,\ldots,x_k)\in \cK_k$, and
$u\in M_{x_1,\ldots,x_k}$.
Then, for all $i\in\{1,\ldots,k\}$, $\lambda_i\in\R^N$ exists so that
\beq 
\label{ML}
I'(u)[\psi]=\int_{B_{R}(x_i)} u^\delta(x)\,
\psi(x)\, (\lambda_i\cdot (x-x_i))\, dx\qquad\forall\psi\in
H_0^1(B_{R}(x_i)).
\eeq
\end{prop}

Next step in the program of building $\forall k\in\N\setminus\{0\}$ a
good candidate critical point emerging around $k$ points is to show
that the supremum of $\mu(x_1,\ldots,x_k)$, when $(x_1,\ldots,x_k)$
varies in $\cK_k$, is achieved.

\vspace{2mm}

For all $k\in\N\setminus\{0\}$ we set 
$$
\mu_k=\sup_{\cK_k}\mu(x_1,\ldots,x_k)=\sup_{\cK_k}\min_{S_{x_1,\ldots,x_k}}
I(u),
$$
$$
\mu^\infty_k=\sup_{\cK_k}\mu^\infty(x_1,\ldots,x_k)=\sup_{\cK_k}\min_{S_{x_1,\ldots,x_k}}
I^\infty(u).
$$
We first consider the case of functions emerging just around one point, then we prove the  wanted result in the case $k > 1.$

We start recalling a lemma, proved in \cite{CPS1} that characterizes the set $M^\infty_y$ and states $\mu_1^\infty$ is achieved. 

\begin{lemma}
\label{4.2CPS}
The relation
\beq
\label{3.5.a}
\mu_1^\infty=m_\infty=\mu^\infty (y)
\eeq
holds for all $y\in\R^N$, moreover
$$
M^\infty_y=\{w_y(x)\}
$$
$w_y$ being the function defined in (\ref{N2'.1}).
\end{lemma}

Next two lemmas, which are basic for the proof of the subsequent
Propositions \ref{GP3.8} and \ref{GP3.9}, describe some feature of
minimizers on classes $S_{x_n,}$ when the sequence of points
$\left\lbrace x_n \right\rbrace$ goes to infinity. First one states
the action of such a sequence of functions converges to $m_\infty,$
the action of the   ground state  $w$ of the limit problem, the second
one asserts that the shape of the functions as $n\to \infty$
approaches the shape of the ground state $w$. 

\begin{lemma}
\label{minfty}
Let $\{x_n\}_n$, $x_n\in\R^N$, be such that $|x_n|\lo +\infty$, and
let $u_n\in M_{x_n}$.
Let $\tilde w_{x_n}$ be the projection of $w_{x_n}$ on $S_{x_n}$ and
$\hat u_n$ the projection on $S^\infty_{x_n}$ of $u_n$.
Then
\beq
\label{a_1}
\mu(x_n)\lo m_\infty;
\eeq
\beq
\label{c_1}
\lim_{n\to  + \infty} I(u_n)=\lim_{n\to  + \infty} I(\hat u_n)=\lim_{n\to  + \infty} I(\tilde w_{x_n})=m_\infty.
\eeq
\end{lemma}

\no\proof Proposition \ref{L7*} implies that $\{\|\tilde
w_{x_n}\|\}_n$ is bounded, then, by using the asymptotic decay of $w$
and assumptions $(h_1),\, (h_2),\, (h_3)$, we deduce
\beq
\label{4.1}
\begin{array}{rcl}
\vspace{2mm}
\mu(x_n)=I(u_n)\le I(\tilde w_{x_n})& = & I_\infty(\tilde w_{x_n})-{1\over
  2}\int_{\R^N}\alpha(x) \tilde w_{x_n}^2dx+{1\over
  q+1}\int_{\R^N}b(x) \tilde w_{x_n}^{q+1}dx\\
& \le &m_\infty +o(1).
\end{array}
\eeq
On the other hand
\beq
\label{4.3}
\begin{array}{rcl}
\vspace{2mm}
\mu(x_n) & = & I(u_n)\ge I(\hat u_n)
\\ \vspace{2mm}
& = & I_\infty(\hat u_n)-{1\over
  2}\int_{\R^N}\alpha(x) \hat u_n^2dx+{1\over
  q+1}\int_{\R^N}b(x) \hat u_n^{q+1}dx
\\
& \ge &m_\infty-{1\over
  2}\int_{\R^N}\alpha(x) \hat u_n^2dx+{1\over
  q+1}\int_{\R^N}b(x) \hat u_n^{q+1}dx.
\end{array}
\eeq
Therefore, if we show that
\beq
\label{4.2}
 \int_{\R^N}\alpha(x) \hat u_n^2dx=o(1),
\qquad
 \int_{\R^N}b(x) \hat u_n^{q+1}dx=o(1)
\eeq
(\ref{a_1}) and (\ref{c_1}) are proved.

Let us fix $\e>0$, then by (\ref{as}) a number $R_\d>R$ can be found
such that $\forall \rho>R_\d$ and for large $n$
\beq
\label{4.6}
\begin{array}{l}
\vspace{2mm}
\int_{\R^N\setminus B_\rho(x_n)}\alpha(x) \hat u_n^2\,dx\le|\alpha|_\infty
\int_{\R^N\setminus B_\rho(x_n)}  u_n^2\,dx<\e,\\
\int_{\R^N\setminus B_\rho(x_n)}b(x) \hat u_n^{q+1}dx\le|b|_\infty
\int_{\R^N\setminus B_\rho(x_n)}  u_n^{q+1}dx<\e.
\end{array}
\eeq
Moreover, thanks to $(h_1),\, (h_2),\, (h_3)$,  and $|x_n|\lo +\infty$, we
can assert that $\bar \rho>R_\d$ exists so that  for large $n$
\beq
\label{4.7}
|\alpha|_{L^\infty(B_{\bar\rho}(x_n))}<{\e},\qquad
|b|_{L^\infty(B_{\bar\rho}(x_n))}<{\e}.
\eeq
So, we have
\beq
\label{4.4}
\begin{array}{rcl}
\vspace{2mm}
\int_{\R^N}\alpha(x) \hat u_n^2\, dx& = &\int_{\R^N\setminus B_{\bar\rho}(x_n)}\alpha(x)  u_n^2\, dx+
\int_{ B_{\bar\rho}(x_n)}\alpha(x) \hat u_n^2\, dx\\
\vspace{3mm}
& <& \e+{\e} \int_{ B_{\bar\rho}(x_n)} \hat
u_n^2\, dx,\\ \vspace{2mm}
\int_{\R^N}b(x) \hat u_n^{q+1}dx& = & \int_{\R^N\setminus B_{\bar\rho}(x_n)}b(x)  u_n^{q+1}dx+
\int_{ B_{\bar\rho}(x_n)}b(x) \hat u_n^{q+1} dx\\
&<&\e+{\e} \int_{ B_{\bar\rho}(x_n)} \hat
u_n^{q+1}dx.
\end{array}
\eeq
Hence (\ref{4.2}) follows if we prove that $\{\|\hat u^\d_n\|\}_n$ is
bounded.
To do this, first we observe that the same argument displayed in the
inequalities (\ref{8.-1}), (\ref{8.0}) of Proposition \ref{Pmin}
together  with $I(u_n)\le m_\infty+o(1)$, allow us to conclude that
$\{\| u_n^\d\|/|u^\d_n|_{p+1}\}_n$ is bounded.
This done, the boundedness of $\{\|u^\d_n\|\}_n$ is easily obtained, because
otherwise the impossible inequality would be true
$$
0<J(u_n^\d)\le c_1\|u_n^\d\|^2-c_2|u_n^\d|_{p+1}^{p+1}+c_3\to -\infty.
$$
Moreover, $I'(u_n)[u_n^\d]=0$ and (\ref{delta}) imply
$$
\begin{array}{rcl}
\vspace{2mm}
0 & \ge & c_1\|u^\d_n\|^2+(a_0\d-2^{p-1}\d^p)\int_{\supp
  u^\d_n}u^\d_n dx-2^{p-1}|u^\d_n|^{p+1}_{p+1} \\
& \ge & c_1\|u_n^\d\|^2-c_2\|u^\d_n\|^{p+1}
\end{array}
$$
from which
\beq
\label{paura}
\|u^\d_n\|\ge c>0 \qquad\forall n\in\N, \quad c\in\R.
\eeq
Lastly, using (\ref{paura}) and arguing as in (\ref{pr}), (\ref{pr2}),
we obtain that $\{\theta^\infty(u_n)\}_n$, and hence $\{\|\hat u_n^\d\|\}_n$,
is bounded. 

\qed

\begin{lemma}
\label{L3.6*}
Let $\{x_n\}_n$, $x_n\in\R^N$ be such that $|x_n|\lo + \infty$.
Let $\{u_n\}_n$ be a sequence such that for all $n$, $u_n\in S_{x_n}$
and
\beq
\label{G*bis}
0<u_n<c\, e^{-\tau\dist (x,\supp u^\d_n)}
\eeq
for some $\tau>0$ and $c >0$  not depending on $n$.
Let $\hat u_n$ be the projection of $u_n$ on $S^\infty_{x_n}$ and let
the relation
\beq
\label{G**}
\lim_{n\to+\infty}I(u_n)=\lim_{n\to+\infty}I(\hat u_n)=m_\infty
\eeq
holds true.
Then
\beq
\label{3.11}
\begin{array}{ccl}
(i) & & u_n(\cdot +x_n)\lo w\qquad\mbox{ in }H^1(\R^N)\\
 (ii) & &\hat u_n(\cdot +x_n)\lo w\qquad\mbox{ in }H^1(\R^N).
\end{array}
\eeq
\end{lemma}

\no\proof
Arguing exactly as in the proof of Lemma \ref{minfty} we deduce that
$\{\|\hat u_n^\d\|\}_n$ is bounded, $\|\hat u_n^\d\|\ge c>0$, $\forall
n\in\N$ and moreover 
\beq
\label{B1bis}
m_\infty\le \lim_{n\to +\infty}I_\infty(\hat u_n)=\lim_{n\to +\infty}I(\hat u_n)=m_\infty.
\eeq
Setting $v_n(x)=\hat u_n(x+x_n)$, $v_n\in S^\infty_0$, and by
(\ref{B1bis}) we deduce
\beq
\label{B3}
\lim_{n\to +\infty}I_\infty(\hat u_n)=\lim_{n\to +\infty}I_\infty (v_n)=m_\infty.
\eeq
Furthermore $\|v_n\|=\|\hat u_n\|$, so  $\{\|v_n\|\}_n$ is bounded, $\| v_n^\d\|\ge \bar{c}>0,$ and 
we can assert the existence of a function $v\in H^1(\R^N)$ such that,
up to a subsequence, 
\beq
\label{B2}
v_n\lo v\qquad\left\{ 
\begin{array}{ll}
 \mbox{ weakly in
}H^1(\R^N) & (i)\\
\mbox{ a.e. in }\R^N & (ii)\\ 
\mbox{ strongly in }L^{s}_{\loc}(\R^N),\ 2\le s < 2^*. & (iii)
\end{array}\right.
\eeq
To prove (\ref{3.11})$(ii)$ we intend to show that $v_n\to v$ strongly
in $H^1(\R^N)$ and $v\in M^\infty_0$, that is $v=w$ by Lemma \ref{4.2CPS}.
We start observing that  (\ref{B2})$(ii)$ gives $v\le \d$ in $\R^N\setminus B_R(0)$ and that
(\ref{B2})$(iii)$ together with assumption
(\ref{G*bis}) implies
\beq
\label{B4}
v_n\to v\quad\mbox{ strongly in }L^s(\R^N),\qquad 2\le s < 2^*.
\eeq
Now, since by (\ref{B1bis}) $\lim_{n\to +\infty} I_\infty (v_n) = m_\infty,$ 
arguing as in (\ref{8.-1}), (\ref{8.0}), with $I_\infty$ and
$J_\infty$ in place of $I$ and $J$, we deduce that the sequence
$\{\|v_n^\d\|/|v_n^\d|_{p+1}\}_n$ is bounded, from which
$|v_n^\d|_{p+1} > \hat{c}> 0$ and, then by (\ref{B4}),
$|v^\d|_{p+1}\neq 0$. 
Again using the  strong convergence of $v_n$ to $v$ in $L^{p+1}(\R^N)$
and the definition of $S_0^\infty$ we get    
$$
m_\infty
\le I_\infty( v_\delta+\theta^\infty( v)v^\delta)
\le \liminf_{n\to  + \infty}I_\infty((v_n)_\delta+\theta^\infty(v)v_n^\d)
$$
\beq
\le\liminf_{n\to  + \infty} I_\infty((v_n)_\d +v_n^\d)
= \lim_{n\to  + \infty} I_\infty (v_n)=m_\infty,
\label{1107}
\eeq
from which 
\beq
\label{70}
\lim_{n\to+\infty}I_\infty((v_n)_\d+\theta^\infty(v)v^\d_n)=I_\infty( v_\delta+\theta^\infty( v)v^\delta)
\eeq
 clearly follows.

(\ref{70}) implies that
 $v_n\to  v$ strongly in $H^1(\R^N)$, so $v \in S_0^\infty$ and
$m_\infty = I_\infty (v)$.  Hence $ v\in M_0^{\infty}$,
 Lemma \ref{4.2CPS} gives $v=w$, and  (\ref{3.11})$(ii)$ is proved.
 
Lastly (\ref{3.11})$(i)$ follows from (\ref{3.11})$(ii)$   arguing as in (\ref{pr}) and (\ref{pr2}) to show that the sequence $\left\{\theta^\infty (u_n)\right\}$ of the coefficients of the projections of $u_n$ on $S^\infty_0$ is bounded and goes to $1$ as $n \to \infty.$

\qed

\begin{prop}
\label{GP3.8}
Let assumption $(h_4)$ be satisfied.
Then
\beq
\label{m1}
\mu_1>m_\infty.
\eeq
Moreover $\mu_1$ is achieved, i.e. 
\beq
\label{G3.25}
\exists\, \bar  x\in\R^N\ \mbox{ such that }\ \mu(\bar x)=\mu_1.
\eeq
\end{prop}

\noindent \proof
Let $\zeta_n=n\,\bar \zeta$ and let
$\{u_n\}_n$ be such that $u_n\in M_{\zeta_n}$.
Set $\hat u_n=(u_n)_\d+\theta^\infty(u_n)u_n^\d$, $\hat u_n\in
S^\infty_{\zeta_n}$.

Since, by Lemma \ref{minfty}, $\mu(\zeta_n)\to m_\infty$ and by definition of $S_{\zeta_n}$
 $\max_{t>0} I((u_n)_\d+tu_n^\d)=I((u_n)_\d+u_n^\d)=I(u_n)$ we deduce 
\beq
\label{16.0}
\begin{array}{rcl}
\vspace{2mm}
\mu_1&\ge& I(u_n)  \ge I(\hat u_n)=I_\infty(\hat u_n)-{1\over
  2}\int_{\R^N}\alpha(x)\hat u_n^2dx+{1\over
  q+1}\int_{\R^N}b(x)\hat u_n^{q+1}dx\\
\vspace{2mm}
& \ge & m_\infty-{1\over
  2}\left[\int_{\cD}\alpha(x)\hat u_n^2dx+\int_{\R^N\setminus\cD}
  \alpha(x)\hat u_n^2dx\right] +{1\over
  q+1}\int_{B_{1}(\zeta_n)}b(x)\hat u_n^{q+1}dx.
\end{array}
\eeq
Now, thanks to $(h_4)$ we get
\beq
\label{16.1b}
\lim_{n\to +\infty} \left[\int_{B_{1}(\zeta_n)}b(x)\hat
  u_n^{q+1}dx\right]\, 
e^{\eta\, n}=+\infty.
\eeq
On the other hand, considering that  $d\ge {1/2}$ and the asymptotic
decay (\ref{as}) we obtain
\beq
\label{16.2b}
\int_{\R^N\setminus\cD} \alpha(x)\hat u_n^2dx\le 
\int_{\R^N\setminus B_{nd}(\zeta_n)} \alpha(x)\hat u_n^2dx\le 
c\, e^{-\eta\, n},
\eeq
while, using  $(h_4)$ and (\ref{as}) we deduce
\beq
\label{16.3}
\begin{array}{rcl}
\vspace{2mm}
\left[\int_{\cD}\alpha(x)\hat u_n^2dx\right]\, e^{\eta\, n} & \le &
c\int_{\cD}e^{-\eta|x|+\eta\, n}\hat u_n^2dx
\\ \vspace{2mm}& = & 
c\int_{\cD-\{\zeta_n\}}e^{\eta(|\zeta_n|-|x+\zeta_n|)}\hat
u_n^2(x+\zeta_n)\, dx
\\ \vspace{2mm}& \le & 
c\int_{\R^N}e^{\eta\, |x|}\hat u_n^2(x+\zeta_n)\, dx
\\ & < & +\infty.
\end{array}
\eeq
Thus, (\ref{m1}) follows directly from (\ref{16.0})--(\ref{16.3}).

To prove (\ref{G3.25}), let us consider a sequence $\{y_n\}_n$,
$y_n\in\R^N$, such that $\lim_{n\to + \infty}\mu(y_n)$ $=\mu_1$ and let $\{v_n\}_n$ be
such that $v_n\in M_{y_n}$.
By Lemma  \ref{minfty}, $\{y_n\}_n$ is bounded, so up to a subsequence
$y_n\to\bar x$ and taking into account $\bar u\in M_{\bar x}$, we have
$I(\bar u)=\mu(\bar x)\le \mu_1$.
On the other hand, let us consider $\bar u_n(x)=\bar u(x-y_n+\bar x)$
and $\bar v_n=(\bar u_n)_\d+\theta(\bar u_n)u^\d_n$, then $\bar v_n\in
S_{y_n}$ and, in view of the continuity of $\theta(\bar u_n(x-y_n+\bar
x))$ with respect to $y_n$, we get 
$$
\mu_1=\lim_{n\to + \infty}I(v_n)\le\lim_{n\to + \infty}I(\bar v_n)=I(\bar u)
$$
so $\mu(\bar x)=I(\bar u)=\mu_1$.

\qed

We turn now to the case $k> 1$. Before proving the desired result in
Proposition \ref{GP3.9} it is useful to show the upper semi-continuity
of the map $( x_1,\ldots,x_k) \to \mu(x_1,\ldots,x_k ).$ 
\begin{lemma}
\label{usc}
Let $\{( x_1^n,\ldots,x^n_k)\}_n$ be a sequence of $k$-tuples
belonging to $\cK_k$  such that 
$$
\lim_{n\to + \infty}( x_1^n,\ldots, x^n_k)=
( x_1,\ldots, x_k).
$$
Then 
\beq
\label{g3.71}
\limsup_{n\to + \infty} \mu(x_1^n,\ldots,x^n_k)\le\mu(x_1,\ldots,x_k ).
\eeq
\end{lemma}

\noindent\proof
Let us consider 
$\bar u\in M_{x_1,\ldots,x_k}$.
Set for all $n\in\N$
$$
u_n=(u_n)_\d+\sum_{i=1}^k\theta_i(u_n)(u_n)_i^\d
$$
where for all $i\in\{1,\ldots,k\}$
$$
(u_n)^\d_i=\bar u^\d_i(x+x_i-x_i^n)
$$
and $(u_n)_\d$ is the unique positive minimizer for the minimization
problem
$$
\min \{I(u)\ :\ u\in H^1(\R^N),\ 0\le u \le \delta,\ u=\delta\mbox{ on
}\cup_{i=1}^k \supp (u_n)^\d_i\}.
$$
Remark that such a minimizer exists and is unique because the same
argument of Lemma \ref{L2.7} shows that $I$ is coercive and convex on
the convex set $\{ u\in H^1(\R^N)\ :\ 0\le u \le \delta,\ u=\delta\mbox{ on
}\cup_{i=1}^k \supp (u_n)^\d_i\}$.

Then we have $u_n\in S_{ x_1^n,\ldots,x^n_k}$, $u_n\to\bar u$ in
$H^1(\R^N)$, and
$$
\limsup_{n\to + \infty}\mu(x_1^n,\ldots,x^n_k )\le 
\limsup_{n\to + \infty}I(u_n)=I(\bar u)=\mu(x_1,\ldots,x_k),
$$
as desired.

\qed

\begin{prop}
\label{GP3.9}
Let assumption  $(h_4)$ be satisfied.
For all $k\in \N\setminus\{0\}$
\beq
\label{G3.91}
\left\{
\begin{array}{rl}
i) & \exists (\bar x_1,\ldots,\bar x_k)\in\cK_k\ :\ \mu_k=\mu( \bar
x_1,\ldots,\bar x_k)\\
ii) & \mu_{k+1}>\mu_k+m_\infty.
\end{array}
\right.
\eeq
\end{prop}

\noindent \proof
The proof is carried out by an inductive argument on $k$.

{\underline {Step 1}}\qquad $k=1$.\quad 
In this case relation (\ref{G3.91}){\em (i)} is just relation
(\ref{G3.25}) stated in Proposition \ref{GP3.8}.
Let us prove that
\beq
\label{G3.27}
\mu_2>\mu_1+m_\infty.
\eeq
For all $n$, let us consider $u_n\in M_{\bar x,\zeta_n}$ where
$\zeta_n=n\bar \zeta$ and $\bar
x\in\R^N$ is such that $\mu_1=\mu(\bar x)$; $u_n$ can be written as
$u_n=(u_n)_\d+\bar u^\d_n+\tilde u^\d_n$ where $\bar u^\d_n$ and
$\tilde u^\d_n$ are the emerging parts around $\bar x$ and $\zeta_n$
respectively. 
We remark that, for large $n$,
$$
B_R(\bar x)\subset \{x\in\R^N\ :\
  (x\cdot{\zeta_n\over|\zeta_n|})<{|\zeta_n|\over
    2}-1\}, 
\ 
B_R(\zeta_n)\subset \{x\in\R^N\ :\
  (x\cdot{\zeta_n\over|\zeta_n|})>{|\zeta_n|\over
    2}+1\}. 
$$
Set now
$$
\Sigma_n=\left\{x\in\R^N\ :\ {|\zeta_n|\over
    2}-1<\left(x\cdot{\zeta_n\over|\zeta_n|}\right)<{|\zeta_n|\over
    2}+1\right\},
$$
\beq
\label{chi}
\begin{array}{c}
\vspace{2mm}
\chi_n(x)
=\chi\left(\left|\left(x\cdot{\zeta_n\over|\zeta_n|}\right)-{|\zeta_n|\over
    2}\right|\right)
\\ 
\chi\in C^\infty(\R,[0,1]) \mbox{ s.t. }
\chi(t)=0\ \mbox{ if }|t|\le 1/2, \ \chi(t)=1\ \mbox{ if }\ |t|\ge
1,
\end{array}
\eeq
define
$$
v_n(x)=\chi_n(x)u_n(x)
$$
and evaluate $I(v_n)$.
Taking into account the exponential decay of $u_n$, we infer
\beq
\label{G3.29}
\begin{array}{rcl}
\vspace{2mm}
I(v_n) & \le & I(u_n)+{1\over
  2}\int_{\Sigma_n}\left[|\D\chi_n|^2-{1\over 2}\Delta
  \chi_n^2\right]u_n^2dx
\\ \vspace{2mm}
& & -{1\over q+1}\int_{\Sigma_n}(1-\chi_n^{q+1})\, b(x)\, u_n^{q+1}dx
 +{1\over p+1}\int_{\Sigma_n}(1-\chi_n^{p+1})u_n^{p+1}dx
\\ \vspace{2mm}
&\le& \mu_2 +c_1\int_{\Sigma_n}(u_n)_\delta^2dx
+c_2\int_{\Sigma_n}(u_n)_\delta^{q+1}dx
+c_3\int_{\Sigma_n}(u_n)_\delta^{p+1}dx
\\ \vspace{2mm}
& \le & \mu_2+o(e^{-\eta|\zeta_n|}).
\end{array}
\eeq

On the other hand 
$$
v_n(x)=v_n^I(x)+v_n^{II}(x)
$$ 
with 
$$
v^I_n(x)\ =\ \left\{ \begin{array}{lc}
0 & \mbox{ if } \left(x\cdot{\zeta_n\over|\zeta_n|}\right)\ge {|\zeta_n|\over
    2}-{1\over 2}\\
\bar u^\d+\chi_n(x)(u_n)_\d(x) 
& \mbox{ if } \left(x\cdot{\zeta_n\over|\zeta_n|}\right)< {|\zeta_n|\over
    2}-{1\over 2}
\end{array}\right.
$$
and
$$
v^{II}_n(x)\ =\ \left\{ \begin{array}{lc}
0 & \mbox{ if } \left(x\cdot{\zeta_n\over|\zeta_n|}\right)\le {|\zeta_n|\over
    2}+{1\over 2}\\
\tilde u^\d+\chi_n(x)(u_n)_\d(x) 
& \mbox{ if } \left(x\cdot{\zeta_n\over|\zeta_n|}\right)>{|\zeta_n|\over
    2}+{1\over 2}.
\end{array}\right.
$$
By definition 
$v_n^I\in S_{\bar x}$, $v_n^{II}\in S_{\zeta_n}$, and
\beq
\label{11.2}
I(v_n^I)\ge\mu(\bar x)=\mu_1.
\eeq
To estimate $I(v^{II}_n)$, let us consider, for all $n$, $z_n\in
M_{\zeta_n}$.
We have 
\beq
\label{G+}
I(v^{II}_n)\ge I(z_n).
\eeq
On the  other hand, taking into account (\ref{as}) and that 
$(\supp v^I_n)\cap (\supp v^{II}_n)=\emptyset$ we get
\beq
\label{natale}
I(v^I_n)+I(v^{II}_n)=I(v_n)=I(u_n)+o(1)
\eeq
while, using the asymptotic decay of $v^I_n$ and $z_n$, we can assert
that for large $n$ $v^I_n+z_n\in S_{\bar x,\zeta_n}$ and that 
$$
I(u_n)\le I(v^I_n+z_n)=I(v^I_n)+I(z_n)+o(1).
$$
So for large $n$ we have
$$
I(v^{II}_n)\le I(z_n)+o(1),
$$
that, together with (\ref{G+}) and (\ref{c_1}) of Lemma \ref{minfty}
implies
$$
m_\infty=\lim_{n\to + \infty} I(z_n)=\lim_{n\to + \infty} I(v^{II}_n)
$$
and, by Lemma \ref{L3.6*}
$$
v_n^{II}(\cdot +\zeta_n)\lo w\mbox{ in } H^1(\R^N)
$$
$$
\hat v_n''(\cdot+\zeta_n)\lo w \mbox{ in } H^1(\R^N)
$$
being $\hat v^{II}_n$ the projection of $v^{II}_n$ on $S^\infty_{\zeta_n}$.
Now, working as in (\ref{16.0})--(\ref{16.3}), we deduce
\beq
\label{newr}
I(v_n^{II})=m_\infty+l_ne^{-\eta\, n}+O(e^{-\eta\, n})\quad\mbox{ with
}l_n\lo +\infty.
\eeq
Collecting (\ref{G3.29}), (\ref{natale}), (\ref{11.2}) and (\ref{newr})
we finally obtain
$$
\begin{array}{rcl}
\vspace{2mm}
\mu_2& \ge&
I(v_n)+o(e^{-\eta |\zeta_n|})\\
&\ge & \mu_1+m_\infty+l_ne^{-\eta\, n}+O(e^{-\eta\, n})
\end{array}
$$
that gives (\ref{G3.27}).

{\underline {Step 2}}\qquad $k>1$.\quad 
We assume (\ref{G3.91}){\em (i)},{\em (ii)} hold true for all $h<k$.
Let us first show that 
\beq
\label{G3.9a} 
\exists (\bar x_1,\bar x_2,\ldots,\bar x_k)\in\cK_k\ :\ \mu_k=\mu(\bar
x_1,\bar x_2,\ldots,\bar x_k).
\eeq
Let $\{(x^n_1,\ldots,x^n_k)\}_n$ be a sequence of elements of $\cK_k$
such that 
\beq
\label{G3.9*}
\lim_{n\to + \infty} \mu(x_1^n,\ldots,x^n_k)=\mu_k,
\eeq
we claim that it is bounded. 
Indeed, once proved the claim, we can assert that, up to a
subsequence, $\{(x^n_1,\ldots,x^n_k)\}_n$ is convergent, moreover,
denoting by $(\bar x_1,\ldots,\bar x_k)$ its limit, we have
$\mu_k=\mu(\bar x_1,\ldots,\bar x_k)$ by Lemma \ref{usc}.
To prove the claim we argue by contradiction and we assume that for
some $j$, $0\le j < k$, $k-j$ sequences among $\{x^n_i\}_n$,
$i\in\{1,\ldots,k\}$, are unbounded. 
Without any loss of generality, we can suppose that, up to
subsequences,
$$
\lim_{n\to + \infty}|x^n_i|=\infty\qquad\forall i\in\{j+1,\ldots,k\}.
$$
Let us fix $s_n\in M_{x^n_1,\ldots,x^n_j}$, consider 
$\tilde w_{x^n_{j+1}}\vee \tilde w_{x^n_{j+2}}\vee \ldots\vee
\tilde w_{x^n_{k}}\in S_{x^n_{j+1},\ldots,x^n_k}$ and define $z_n\in 
S_{x^n_{1},\ldots,x^n_k}$ as
$$
z_n=(z_n)_\d+\sum_{i=1}^j(s_n)^\d_i+(\tilde w_{x^n_{j+1}}\vee \ldots\vee
\tilde w_{x^n_{k}})^\d
$$
with
$$
(z_n)_\d=\left[1-\chi\left({2\over
      d_n}|x|\right)\right](s_n)_\d+\chi\left({2\over
    d_n}|x|\right)(\tilde w_{x^n_{j+1}}\vee\cdots\vee\tilde w_{x^n_k})_\d,
$$
where $d_n=\min\{|x^n_{i}|\ :\ i\in\{j+1,\ldots,k\}\}$ and $\chi$ is
the cut-off function introduced in (\ref{chi}).
For large $n$, we have
\beq
\label{Gpallino}
\begin{array}{rcl}
\vspace{2mm}
\mu(x^n_1,\ldots,x^n_k)&\le& I(z_n)=I((z_n)_\d)+J(z_n^\d)
\\ \vspace{2mm}
& = &  I((z_n)_\d)+J(s_n^\d)+J\left((\tilde w_{x^n_{j+1}}\vee\cdots\vee \tilde
w_{x^n_k})^\d\right)\\ \vspace{2mm}
&=& I(s_n)+I(\tilde w_{x^n_{j+1}}\vee\cdots\vee \tilde
w_{x^n_k})+I((z_n)_\d)-I((s_n)_\d)\\
& & -I\left((\tilde w_{x^n_{j+1}}\vee\cdots\vee \tilde
w_{x^n_k})_\d\right).
\end{array}
\eeq 
Now
\beq
\label{G+bis}
I(s^\d_n)\le\mu_j,
\eeq
moreover, considering that $\tilde w_{x^n_{j+1}}\wedge (\tilde
w_{x^n_{j+2}}\vee\ldots \vee\tilde w_{x^n_{k}})\le\d$ and hence by Lemma
\ref{L2.7} $I(\tilde w_{x^n_{j+1}}\wedge (\tilde
w_{x^n_{j+2}}\vee\ldots \vee\tilde w_{x^n_{k}}) )>0$
$$
\begin{array}{rl}
\vspace{2mm}
\hspace{-1mm} I(\tilde w_{x^n_{j+1}}\vee\ldots\vee \tilde w_{x^n_{k}})
& \hspace{-2mm} = \hspace{-1mm} I(\tilde w_{x^n_{j+1}})\hspace{-1mm}+\hspace{-1mm}
I(\tilde
w_{x^n_{j+2}}\vee\ldots \vee\tilde w_{x^n_{k}} )
 \hspace{-1mm}- \hspace{-1mm}I(\tilde
w_{x^n_{j+1}}\wedge 
(\tilde w_{x^n_{j+2}}\vee\ldots\vee \tilde w_{x^n_{k}} ))\\
&\hspace{-2mm} \le  I(\tilde w_{x^n_{j+1}})+ I(\tilde w_{x^n_{j+2}}\vee\ldots \vee\tilde w_{x^n_{k}} )
\end{array}
$$
from which, repeating the argument on
$I(\tilde w_{x^n_{j+2}}\vee\ldots \vee \tilde w_{x^n_{k}})$, 
$I(\tilde w_{x^n_{j+3}}\vee\ldots \vee\tilde w_{x^n_{k}})$, \ldots,
$I(\tilde w_{x^n_{k-1}}\vee \tilde w_{x^n_{k}})$
we obtain, using Lemma \ref{minfty},
\beq
\label{G*}
I(\tilde w_{x^n_{j+1}}\vee\ldots \vee \tilde w_{x^n_{k}} )
\le\sum_{i=j+1}^k(\tilde w_{x^n_i})= (k-j)\, m_\infty+o(1).
\eeq
Lastly the definition and the asymptotic decay of $z_n$, $s_n$, 
$\tilde w_{x^n_{j+1}}\vee\ldots \vee \tilde w_{x^n_{k}}$ allow to
conclude 
\beq
\label{G.}
I((z_n)_\d)-I((s_n)_\d)-
I(\tilde w_{x^n_{j+1}}\vee\ldots \vee \tilde w_{x^n_{k}} )=o(1).
\eeq
 Inserting (\ref{G+bis}), (\ref{G*}) and (\ref{G.}) in
 (\ref{Gpallino}) we get
$$
\mu(x^n_1,\ldots,x^n_k)\le \mu_j+(k-j)m_\infty+o(1)
$$
that, with (\ref{G3.9*}) implies
$$
\mu_k\le \mu_j+(k-j)m_\infty
$$
contradicting $\mu_k>\mu_j+(k-j)m_\infty$, that comes from
(\ref{G3.91}){\em (ii)}.
Thus the claim is proved.

To conclude the proof, relation (\ref{G3.91})$(ii)$ is left to be shown.
Since the argument is similar to that of step 1, we skip some details.

Let us consider $u_n\in M_{\bar x_1,\ldots, \bar x_k,\zeta_n}$ where
$(\bar x_1,\ldots, \bar x_k)$ is the above found $k$-tuple for which
$\mu_k=\mu(\bar x_1,\ldots, \bar x_k )$.
$u_n$ can be written as
$$
u_n=(u_n)_\d+\sum_{i=1}^k(u^\d_n)_i+\tilde u_n^\d
$$
where $(u^\d_n)_i$ are the emerging parts around $\bar x_i$,
$i=1,\ldots,k$, and $\tilde u_n^\d$ is the emerging part around
$\zeta_n$.
Remark that, for large $n$, we can assume
$$
\bigcup\limits_{i=1}^k\supp (u^\d_n)_i\subset\bigcup\limits_{i=1}^k
B_R(\bar x_i)\subset 
\left\{x\in\R^N\ :\ \left(x\cdot{\zeta_n\over|\zeta_n|}\right)<{|\zeta_n|\over
    2}-1\right\}
$$
$$
\supp\tilde u^\d_n\subset B_R(\zeta_n)\subset
\left\{x\in\R^N\ :\ \left(x\cdot{\zeta_n\over|\zeta_n|}\right)>{|\zeta_n|\over
    2}+1\right\}.
$$
Setting $v_n=\chi_n(x)u_n(x)$, where $\chi_n$ are the cut-off
functions defined in (\ref{chi}), repeating computations made to
prove (\ref{G3.29}) we obtain
\beq
\label{G3.37}
I(v_n)\le\mu_{k+1}+o(e^{-\eta|\zeta_n|}).
\eeq
On the other hand, as in step 1, we can write
$$
v_n(x)=v^I_n(x)+v^{II}_n(x)
$$
with
$$
v^I_n\in S_{\bar x_1,\ldots,\bar x_k},\qquad v^{II}_n\in S_{\zeta_n}
$$
and arguments quite analogous to those used to prove (\ref{11.2}) and
(\ref{newr}) give respectively
$$
I(v^I_n)\ge \mu_k
$$
$$
I(v^{II}_n)\ge m_\infty+l_ne^{-\eta\, n}+O(e^{-\eta\, n})\quad\mbox{ with
}l_n\lo +\infty,
$$
that together with (\ref{G3.37})  allow to obtain
(\ref{G3.91})$(ii)$.

\qed

\sezione{Behaviour of the max-min as $|b|_\infty\to 0$}

This section contains some basic estimate about the asymptotic
properties of the functions realizing the max-min values of the action
functional as $|b|_\infty$ goes to $0.$ 

In what follows we assume that, besides $(h_1)$, $(h_2)$, $(h_3)$, also
$(h_4)$ holds true.   
Moreover we consider the family 
$$
\cF=\left\{b\in L^\infty(\R^N)\ :\ b\ge 0,\ |b|_\infty<B_1 \right\}
$$
and we work on problems $({P})$ with variable $b$.
Thus, if $b_n\in\cF$ we denote by
$$
I_n(u)={1\over 2}\int_{\R^N}(|\D u|^2dx+a(x)u^2)dx+{1\over
  q+1}\int_{\R^N}b_n(x)\, |u|^{q+1}dx
-{1\over p+1}\int_{\R^N}|u|^{p+1}dx
$$
and by $J_n$, $S^n_{x_1,\ldots,x_k}$, $M^n_{x_1,\ldots,x_k}$,
$\mu^n(x_1,\ldots,x_k)$, $\mu^n_k$, $\theta^n$ and $\theta^n_i$ the
objects defined with respect to $I_n$, in the same way 
 $J$, $S_{x_1,\ldots,x_k}$, $M_{x_1,\ldots,x_k}$,
$\mu(x_1,\ldots,x_k)$, $\mu_k$, $\theta$ and $\theta_i$ have been
defined with respect to $I$.

 The following Proposition \ref{P4.1}, the subsequent Corollary
 \ref{GC4.1bis}, and Proposition \ref{PN} analyze how the position of
 the emerging parts ,  the ``bumps'' of the functions realizing the
 max-min, varies as $|b|_\infty$ goes to $0,$ while in Lemmas
 \ref{GL4.4}, \ref{GL4.3} and Proposition \ref{w} the asymptotic shape
 of the same functions is described. 
 
\begin{prop}
\label{P4.1}
To any $r>R>0$ there corresponds a constant $\cB_r\in(0,B_1)$ such
that for all $b\in\cF$, $|b|_\infty<\cB_r$ and for all $k\in\N$, $k\ge 2$,
$$
(x_1,\ldots,x_k)\in\cK_k \quad\mbox{ with }\quad
\mu_k=\mu(x_1,\ldots,x_k)
$$
implies
$$
\min\{|x_i - x_j|\ :\  i \neq j ,\   i,j\in\{1,2,...,k\}\} > r.
$$
\end{prop}

\noindent \proof
We argue by contradiction, so we assume the existence of $r\ge R$ and,
for all $n\in\N$, $b_n\in\cF$, $|b_n|_\infty<1/n$, $k_n\in\N$, 
$(x_1^n,\ldots,x_{k_n}^n)\in \cK_{k_n}$, $u_{n}\in
M_{x_1^n,\ldots,x_{k_n}^n}^n$ such that 
$$
I(u_{n})=\mu_{k_n},\qquad \min\{|x^n_i-x^n_j|\ :\ i\neq j,\ 
i,j\in\{1,\ldots,k_n\}\}\le r.
$$
Whitout any loss of generality we can suppose
\beq
\label{15.1}
|x_1^n-x_2^n|\le r\qquad
\forall n\in\N.
\eeq
To prove the statement we intend to show that (\ref{15.1}) implies
\beq
\label{15.2}
\limsup_{n\to
   + \infty}(\mu^n(x_1^n,\ldots,x_{k_n}^n)-\mu^n(x_2^n,\ldots,x_{k_n}^n))<m_\infty.
\eeq
Indeed, this done, the conclusion easily follows because being
$\mu^n(x_1^n,\ldots,x_{k_n}^n)=\mu_{k_n}^n$ and 
$\mu^n(x_2^n,\ldots,x_{k_n}^n)\le \mu_{k_n-1}^n$, (\ref{15.2}) would
imply, for large $n$,  $\mu_{k_n}^n\le\mu_{k_n-1}^n+m_\infty$
contradicting the opposite strict inequality
$\mu_{k_n}^n>\mu_{k_n-1}^n+m_\infty$ which is true by Proposition
\ref{GP3.9}, whose assumptions are verified for all $b\in\cF$.

Let us consider, for all $n$, $\tilde w_{x_1^n}\vee
v_n$ where $v_n\in M^n_{x_2^n,\ldots,x_{k_n}^n}$ and 
\beq
\label{Gruota}
\tilde w_{x^n_1}=(w_{x^n_1})_\d+\theta^n(w_{x^n_1})w_{x^n_1}^\d.
\eeq
Clearly, $\tilde w_{x_1^n}\vee v_n\in S^n_{x_1^n,\ldots,x_{k_n}^n}$, thus
$$
I_n(u_{n})\le I_n(\tilde w_{x_1^n}\vee v_n)=I_n(\tilde
w_{x_1^n})+I_n(v_n)-I_n(\tilde w_{x_1^n}\wedge v_n)
$$
from which
\beq
\label{16.1}
\mu^n(x_1^n,\ldots,x_{k_n}^n)-\mu^n(x_2^n,\ldots,x_{k_n}^n)=I_n(u_{n})-I_n(v_n)\le I_n(\tilde w_{x_1^n})-I_n(\tilde w_{x_1^n}\wedge v_n)
\eeq
follows.

Now, taking into account that
$m_\infty=I^\infty(w)=I^\infty(w_{x^n_1})=\max_{t>0}I^\infty((w_{x^n_1})_\d+t\,
w^\d_{x^n_1})$ and that by Lemma \ref{L7*}
$\{\theta^n(w_{x^n_1})\}_n$ is bounded we have
\beq
\label{16.2}
\begin{array}{rcl}\vspace{2mm}
I_n(\tilde w_{x_1^n})& =& I_\infty(\tilde w_{x_1^n})-{1\over
  2}\int_{\R^N}\alpha(x) \tilde w_{x_1^n}^2dx+{1\over q+1}
\int_{\R^N}b_n(x) \tilde w_{x_1^n}^{q+1}dx\\
&\le &m_\infty+ c |b_n|_\infty\int_{\R^N}w^{q+1}dx+o(1) =  m_\infty+o(1).
\end{array}
\eeq

Moreover, since for all $n$, $I_n$ is easily seen (as in Lemma
\ref{L2.7}) to be coercive on submerged parts and $\tilde
w_{x_1^n}\wedge v_n\le\d$
$$
I_n(\tilde w_{x_1^n}\wedge v_n)=I_n(w_{x_1^n}\wedge v_n)\ge
c\int_{\R^N} (w_{x_1^n}\wedge v_n)^2dx.
$$
Thus, considering that $|x^n_1-x^n_2|<r$, $r>R$ implies 
$\supp(v_n^\d)_2\subset B_{R}(x_2^n)\subset B_{2r}(x_1^n)$, we deduce
\beq
\label{G4.5}
I_n(\tilde w_{x_1^n}\wedge v_n)\ge c \int_{\supp
  (v^\d_n)_2}w^2_{x^n_1}dx\ge c\left(\inf_{B_{2r}(0)}w^2\right)|\supp(v_n^\d)_2|.
\eeq
Therefore, in view of (\ref{16.1}), (\ref{16.2}), and (\ref{G4.5}),
(\ref{15.2}) follows if we show the inequality
\beq
\label{G4.*}
\liminf_{n\to +\infty} |\supp(v_n^\d)_2|>0.
\eeq
Suppose (\ref{G4.*}) false, that is, up to a subsequence, 
\beq
\label{G4.**}
\lim_{n\to +\infty} |\supp(v_n^\d)_2|=0.
\eeq
Then $\{\|(v_n^\d)_2\|/|(v_n^\d)_2|_{p+1}\}_n$ has to be
unbounded; otherwise, in fact, $(v_n^\d)_2(\cdot
-x_2^n)/$ $|(v_n^\d)_2|_{p+1}$ would converge, up to a subsequence, to
some $\tilde u\in H^1_0(B_{R}(0))$  strongly in  $L^{p+1}(B_{R}(0) )$
and a.e. in $B_{R}(0)$. 
Hence $|\tilde u|_{p+1}=1$, while $\tilde u=0$ a.e. in $B_{R}(0)$ by
(\ref{G4.**}) which is impossible.

Therefore, up to a subsequence,  $\|(v_n^\d)_2\|/|(v_n^\d)_2|_{p+1}\to
+\infty$ and, arguing as in (\ref{8.-1}), we get
\beq
\label{G4.6}
\lim_{n\to+\infty} J_n((v_n^\d)_2 )=
\lim_{n\to+\infty} \max_{t>0}J_n(t(v_n^\d)_2 )\ge
\lim_{n\to+\infty} J_n((v_n^\d)_2/|(v_n^\d)_2|_{p+1})=+\infty.
\eeq
Now, set $\tilde v_n=[(v_n)_\d+\sum_{i=3}^{k_n}(v_n^\d)_i]\vee
\tilde w_{x_2^n}$, where $\tilde w_{x^n_2}= (w_{x^n_2})_\d+\theta^n(w_{x^n_2} )
(w_{x^n_2} )^\d\in S^n_{x^n_2}$; then $\tilde v_n\in
S^n_{x_2^n,\ldots,x_{k_n}^n}$ and
\beq
\label{G4.7}
I_n(v_n)\le
I_n(\tilde v_n).
\eeq
On the other hand, analogously to (\ref{16.2}), $I_n(\tilde w_{x_2^n}
)\le m_\infty+o(1)$ thus, considering $J_n((v^\d_n)_2)>0$, we deduce
for large $n$ the inequality
\begin{eqnarray*}
\vspace{2mm}
I_n(\tilde v_n)-I_n(v_n) &\le & I_n((v_n)_\d)+\sum_{i=3}^{k_n}J_n((v_n^\d)_i)+
I_n(\tilde w_{x_2^n})
-[I_n((v_n)_\d)+\sum_{i=2}^{k_n}J_n((v_n^\d)_i)]\\
& = & I_n(\tilde w_{x_2^n} )-J_n((v_n^\d)_2 )
 < 0,
\end{eqnarray*}
that contradicts (\ref{G4.7}) proving (\ref{G4.*}) and completing the
proof.

\qed

\begin{cor}
\label{GC4.1bis}
Let $\{b_n\}_{n}$, $b_n\in\cF$ such that
$\lim_{n\to + \infty}|b_n|_\infty=0$.
Let, for all $n$, $k_n\in\N\setminus\{0\}$,
$(x_1^n,\ldots,x_{k_n}^n)\in\cK_{k_n}$ and $u_n\in H^1(\R^N)$ be such
that 
$$
I_n(u_n)=\mu^n_{k_n},\quad u_n\in M^n_{ x_1^n,\ldots,x_{k_n}^n}.
$$
Then
\beq
\label{G4.1**}
\lim_{n\to + \infty}(\min\{|x_i^n-x^n_j|\ :\ i\neq j,\
i,j\in\{1,\ldots,k_n\}\})=+\infty. 
\eeq
\end{cor}

Arguing as for Proposition \ref{P3.3} the following lemma can be proved.

\begin{lemma}
\label{GL4.1ter}
Let $\{b_n\}_n$ as in Corollary \ref{GC4.1bis}.
Let $k_n\in\N\setminus\{0\}$, $(x_1^n,\ldots,x_{k_n}^n)\in\cK_{k_n}$ and $u_{k_n}
\in M^n_{ x_1^n,\ldots,x_{k_n}^n}$.
Then, setting $d(x)=\dist (x,\supp u^\d_{k_n})$ and fixing $\eta_s\in
(\eta,\sqrt{a_0-\d^{p-1}})$, the relation
\beq
\label{4.12}
0<(u_{k_n})_\d<c\,\d\, e^{-\eta_s d(x)}
\eeq
holds, with $C$ depending neither on $n$ nor on $k_n$.
\end{lemma}

\begin{prop}
\label{PN}
Let $\{b_n\}_{n}$, $b_n\in\cF$, be such that
$\lim_{n\to + \infty}|b_n|_\infty=0$.
Let $\alpha(x)\not\equiv 0$, $\{k_n\}_{n}$, $k_n\in\N\setminus\{0\}$, and
$(x_1^n,\ldots,x_{k_n}^n)\in\cK_{k_n}$ be such that
$
\mu^n_{k_n}=\mu^n(x_1^n,\ldots,x_{k_n}^n)$.

Then
\beq
\label{G4.8}
\lim_{n\to + \infty}\min\{|x_i^n|\ :\ i\in\{1,\ldots,k_n\}\}=+\infty.
\eeq
\end{prop}

\noindent \proof
We prove (\ref{G4.8}) by contradiction.
Without any loss of generality, we assume that $x_1^n\lo
\bar x$, as $n\to +\infty$.

Let us consider $u_n\in M^n_{x_1^n,\ldots,x_{k_n}^n}$, $v_n\in
M^n_{x_2^n,\ldots,x_{k_n}^n}$ and $z_n=v_n\vee \tilde w_{x_1^n}$, $\tilde
w_{x^n_1}$ being defined as in (\ref{Gruota}).
We have
$$
I_n(v_n)=\mu^n(x_2^n,\ldots,x^n_{k_n})\le\mu_{k_n-1}
$$
and, since $|b_n|_\infty\to 0$ and $\alpha\not\equiv 0$
\begin{eqnarray}
\nonumber
\vspace{2mm}
I_n(\tilde  w_{x_1^n}) & \le & I_\infty(\tilde  w_{x_1^n})-{1\over
  2}\int_{\R^N}\alpha(x)\, \tilde  w_{x_1^n}^2dx+o(1)
\\ \vspace{2mm}\label{18.19}
& \leq & I_\infty(w_{\bar x})-{c_1\over 2}\int_{\R^N}\alpha(x) (w_{\bar x})_\d^2dx+o(1)
\\
\nonumber
& < & m_\infty -c_2+o(1),
\end{eqnarray}
with $c_2>0$.
Since  $I_n(v_n\wedge \tilde
w_{x^n_1})>0$, thanks to Lemma \ref{L2.7}, we find
$$
\mu_{k_n}=I_n(u_n)\le I_n(z_n)=I_n(v_n)+I_n(\tilde w_{x^n_1})-I_n( v_n\wedge \tilde w_{x^n_1})< \mu_{k_n-1}+m_\infty-c_2+o(1)
$$
which contradicts (\ref{G3.91}).

\qed

The proof of the following lemmas is analogous to the proof of Lemma
\ref{L3.6*} and \ref{minfty}, respectively.

\begin{lemma}
\label{GL4.4}
Let $\{b_n\}_n$ be as in Proposition \ref{PN}.
Let $\{x_n\}_n,$ $x_n \in \R^N$, be such that $|x_n|\lo  +\infty$ if
$\alpha(x)\not\equiv 0$, and let $u_n\in S^n_{x_n}$ be such that
(\ref{as}) holds with $c$ not depending on $n$.
If 
$$
\lim_{n\to + \infty}I_n(u_n)=m_\infty
$$
then 
$$
u_n(\cdot +x_n)\lo w\qquad\mbox{ in } H^1(\R^N).
$$
\end{lemma}

\begin{lemma}
\label{GL4.3}
Let $\{b_n\}_n$ be as in Proposition \ref{PN}.
Let $\{x_n\}_n,$ $x_n \in \R^N$, be such that $|x_n|\lo  +\infty$ if
$\alpha(x)\not\equiv 0$, and let $u_n\in M^n_{x_n}$.
Let $\tilde w_{x_n}$ be the projection of $w_{x_n}$ on $S^n_{x_n}$ and
$\hat u_n$ the projection on $S^\infty_{x_n}$ of $u_n$.
Then
$$
\mu^n(x_n)\lo m_\infty
$$
and
$$
\lim_{n\to +\infty} I_n(u_n)=\lim_{n\to + \infty} I_n(\hat u_n)
=\lim_{n\to +\infty} I_n(\tilde w_{x_n})=m_\infty.
$$
\end{lemma}

\begin{prop}
\label{w}
Let $\{b_n\}_{n}$ be as in Proposition \ref{PN}.
Let $\{k_n\}_{n}$, $k_n\in\N\setminus\{0\}$.
Let $\{(x_1^n,\ldots,x^n_{k_n})\}_n$ be a sequence of $k_n$-tuples
belonging to $\cK_{k_n}$  and $u_n\in M^n_{x^n_1,\ldots,x^n_{k_n}}$.
Assume that (\ref{G4.1**})  holds true and that 
(\ref{G4.8}) is verified if $\alpha(x)\not\equiv 0$.

Then for all $r>0$
\beq
\label{19.19}
\lim_{n\to +\infty}\sup\{|u_{n}(x+x_i^n)-w(x)|\ :\
i\in\{1,\ldots,k_n\},\ |x|< r\}=0
\eeq
and, for large $n$, 
\beq
\label{19.20}
\supp(u_n^\d)_i\subset \joinrel\subset
B_{R}(x_i^n),\quad \forall i\in\{1,\ldots,k_n\}.
\eeq
\end{prop}

\noindent\proof Relation (\ref{19.20}) is a direct consequence of
(\ref{19.19}), so let us prove (\ref{19.19}).

 First, assume $k_n\ge 2$ $\forall n\in\N$.
Without any loss of generality, in what follows we fix $i=1$.

We consider $\xi_n(x)u_n(x)$, where $\xi_n$ are cut-off functions
defined as $\xi_n(x)=\chi(|x-x_1^n|$ $-{d_n/2})$ with
$\chi$ as in (\ref{chi}) and 
$
d_n= \min\{|x_i^n-x_1^n|\ :
\ i\in\{2,\ldots,k_n\}\}.
$
Since by (\ref{G4.1**})  $d_n\to+\infty$, it is
straightly  verified that, for large $n$
$$
\xi_n(x)u_n(x)=\bar u_n(x)+\check u_n(x)
$$
with $\bar u_n\in S^n_{x_1^n}$, $\check u_n\in
  S^n_{x^n_2,\ldots,x^n_{k_n}}$, $\supp \bar u_n\subset
  B_{{d_n\over 2}-{1\over 2}}(x^n_1)$, and $(\supp \bar u_n)\cap(\supp \check
  u_n)=\emptyset$.

Our argument is carried out by proving the following points
$$
\begin{array}{rl}
\vspace{2mm}
A) & \lim\limits_{n\to +\infty} I_n(\bar u_n)=m_\infty,\qquad \bar u_n(\cdot +
x_1^n)\lo w\ \mbox{ in }H^1(\R^N);\\
B) &  \bar u_n(\cdot + x_1^n)\lo w\ \mbox{ uniformly in } K, \  \forall K
\subset \joinrel\subset B_R(0)\mbox{ compact}.
\end{array}
$$
Actually, once realized $(A)$ and $(B)$, it is not difficult to
conclude.
Indeed the choice of $\d$ and $R$, which implies $w(x)<\d$ when
$|x|>R/2$, the exponential decay (\ref{4.12}), the
relations $u_{k_n}=\bar u_n$ in $B_{d_n/2-1}(x^n_1)$ and
$d_n\to+\infty$, allow us to state that $\bar u_n(\cdot +x^n_1)$
verifies for all $\rho>R$ and for large $n$
$$
-\Delta u+a(\cdot +x_1^n)u+b_n(\cdot +x_1^n)u^q=u^{p}\quad\mbox{ and }\quad 0<u\le c\,
e^{-\eta  |x|}\quad\mbox{ in }B_\rho(0)\setminus B_{R/2}(0)
$$
with a constant $c$ independent of $\rho$.

Hence, taking into account point $(A)$, Lemma \ref{N2}, and 
$a(\cdot +x_1^n)\to a_\infty$, $b_n(\cdot +x_1^n)\to 0$ uniformly in $
B_\rho(0)\setminus B_{R/2}(0),$ regularity arguments give, for all
$\rho>R$, 
$$
 u_n(\cdot + x_1^n)\lo w\quad\mbox{ uniformly in }B_\rho(0)\setminus
B_{R/2}(0)
$$
which, together with point $(B)$ yields (\ref{19.19}).

\noindent {\underline {Proof of point {\em{(A)}}}}\quad Taking into account
Lemma \ref{GL4.1ter}, by computation analogous to those in (\ref{G3.29})
we get 
\beq
\label{G4.33*}
I_n(\bar u_n)+I_n(\check u_n)=I_n(\xi_nu_n)\le I_n(u_n)+O(e^{-\eta
  d_n}).
\eeq
On the other hand, taking $v_n\in M^n_{x^n_1}$, by using again Lemma
\ref{GL4.1ter}, we deduce  for large $n$ the inequality 
$$
I_n(u_n)\le I_n(v_n)+I_n(\check u_n)+o(1)
$$
which together with (\ref{G4.33*}) and $I_n(v_n)\le I_n(\bar u_n)$
gives
$$
\lim_{n\to+\infty}I_n(\bar u_n)=\lim_{n\to+\infty}I_n(v_n).
$$
So, in view of (\ref{G4.8}) and Lemma \ref{GL4.4} 
$$
\lim_{n\to+\infty}I_n(\bar u_n)=m_\infty \ ; \hspace{1,3cm}
u_n(\cdot +x^n_1)\lo w \ \ \mbox{ in }\, H^1(\R^N).
$$

\noindent {\underline {Proof of point {\em{(B)}}}}\quad 
Since $u_n\in M^n_{x^n_1,\ldots,x^n_{k_n}}$ we can apply Proposition
\ref{CPS3.6} that clearly holds for $I_n$, whatever $n\in\N$ is.
Taking into account that $u_n=\bar u_n$ in $B_{R}(x^n_1)$, we can
assert for all $n$ the existence of $\lambda_n\in \R^N$ such that 
\beq
\label{57.11}
I'_n(\bar u_n)[\psi]=\int_{B_{R}(x_1^n)}
\bar u_n^\d(x)\, \psi(x)\, (\lambda_n\cdot (x-x_1^n))\, dx\qquad \forall
\psi\in H^1_0(B_{R}(x_1^n)).
\eeq
Then, a standard bootstrap argument (see e.g. \cite{Ber}) allows us to
state that $\bar u_n(\cdot +x_1^n)\in C^{1,\sigma}(K)$ for all
compact sets $K\subset B_{R}(0)$. 
Moreover, we claim that
\beq
\label{57.12}
\lim_{n\to +\infty}\lambda_n=0.
\eeq
Otherwise, we could assume $|\lambda_n|\ge c>0$ $\forall n\in\N$ and 
$\lim_{n\to+\infty} \lambda_n/ | \lambda_n |=e$, with $|e|=1$.
Testing (\ref{57.11}) with 
$$
\psi_n(x-x_1^n):=  \left({\l_n\over |\l_n|} \cdot (x-x_1^n)\right)
\phi(x-x_1^n),
$$
where $\phi \in C_0^\infty (B_{R}(0))$ is a radially symmetric
function, so that $\phi(x)> 0$ in $B_{R}(0)$, and considering $\bar
u_n(\cdot+x_1^n)\lo w$ in  $H^1(\R^N)$, we get 
$$
\int_{B_{R} (0)}
[(\nabla \bar u_n(x+x_1^n) \cdot \nabla \psi_n (x))
+ a(x+x_1^n)\bar u_n (x+x_1^n)\psi_n (x)]dx
$$
$$
+ \int_{B_{R}  (0)}b_n(x+x_1^n)
(\bar u_n(x+x_1^n))^{q}\psi_n(x)\,dx- \int_{B_{R}  (0)}
(\bar u_n(x+x_1^n))^{p}\psi_n(x)\,dx
$$
$$
 = |\l_n| \int_{B_{R} (0)}
\bar u_n^\d(x+x_1^n)\phi(x)\left(\frac{\l_n}{|\l_n|}\cdot x \right)^2 \,
 dx
$$
$$
=\left(\int_{B_{R}(0)}(w(x))^\d\phi(x)\,(e\cdot x)^2dx+o(1)\right)|\l_n|
$$
\beq
\label{5.17}
\ge c|\l_n|,
\eeq
where $c> 0$.
This is impossible because the first member in the equalities (\ref{5.17})
goes to  $I'_\infty(w)[ (e\cdot x)\, \phi]=0$ as $n\to+\infty$.
Thus (\ref{57.12}) is proved and, by (\ref{57.11}) we can assert that
the sequence $\bar u_n(x+x^n_1)$ is bounded in $\cC^{1,\sigma}(K)$ for
all compact sets $K\subset B_R(0)$ and, then, uniformly converges to
$w$ in all compact sets $K\subset B_R(0)$.

\vspace{2mm}

To conclude the proof we observe that, when  $k_n= 1$ for all  $n\in\N,$ we have just to work as for the points (A) and (B)
with $u_n$ instead of $\bar u_n$. 
Then (A) is nothing but Lemmas \ref{minfty} and
\ref{L3.6*}, while the proof of (B) can be carried out exactly as in the case
$k_n\ge 2$.

\qed

\begin{cor}
\label{C1}
Let $b_n,\ k_n, \ (x_1^n,\ldots,x^n_{k_n})$ and $u_n\in
M^n_{x^n_1,\ldots,x^n_{k_n}}$ be as in Proposition \ref{w}. Then, for
large $n$, the following relation holds true 
\beq
\label{50}
-\Delta
u_n+a(x)u_n+b_n(x)u_n^q=u_n^p+\sum_{i=1}^{k_n}(u_n)_i^\d(\lambda_i^n\cdot
(x-x_i^n)), \ \ \ \ x\in \R^N.
\eeq
for some $\lambda_i^n\in\R^N$,
$i\in\{1,\ldots,k_n\}$.
\end{cor}
\noindent\proof The choice of $R$ and (\ref{19.19}) imply that for
large $n,$ for all $i\in\{1,\ldots,k_n\},$ $\supp(u_n^\d)_i\subset  
B_{R/2}(x_i^n),$ hence by Proposition \ref{P3.3}, $u_n$ is a solution of $-\Delta
u+a(x)u+b_n(x)u^q=u^p$ in $\R^N \setminus
\bigcup_1^{k_n}\bar{B_{R/2}(x_i^n)}.$ On the other hand $u_n$ verifies
(\ref{ML}) in $B_{R}(x_i^n),$ so, by standard regularity results,
(\ref{50})  holds in $B_{\frac{3}{4}R}(x_i^n).$

\sezione{Proof of Theorem \ref{T}}

We have to show that for small $|b|_\infty$ the Lagrange multipliers
$\l_i$ in (\ref{50}) are zero.

Arguing by contradiction, we assume that for all $n\in\N$ there exist
$b_n\in\cF$, $k_n\in\N\setminus\{0\}$, $(x_1,\ldots,x_{k_n})\in
\cK_{k_n}$, and a function $u_n$ such that
\beq
\label{54.1}
|b_n|_\infty\le {1\over n},\quad \mu_{k_n}^n=\mu(x_1,\ldots,x_{k_n})=I_n(u_n),\quad
I'_n(u_n)\neq 0.
\eeq
According to Corollary \ref{C1}, let $\l_i^n$ be the Lagrange multipliers
related to $u_n$, then last inequality in (\ref{54.1}) can be written
$\max\{|\lambda_i^n|$ : $i\in\{1,\ldots,k_n\}\}>0$.

Up to subsequences, we can assume $k_n\equiv \bar k\in\N$ or
$k_n\nearrow + \infty$, $|\lambda_1^n| =\max\{|\lambda_i^n|$ : $i\in\{1,\ldots,k_n\}\}$, $\lim_{n\to
   + \infty}{\lambda_1^n/|\lambda_1^n |}=\bar\lambda$.
Now, let us fix a sequence of real values $\sigma_n>0$ such that
\beq
\label{54.2}
\sigma_n\to 0,\quad \sigma_n k_n\to 0,\quad {\sigma_n \over
  |\l_1^n|}\, k_n \to 0
\eeq
and set
$$
y_1^n=x_1^n+\sigma_n\bar \l ,\qquad y_i^n=x_i^n\quad \forall
i\in\{2,\ldots,k_n\}.
$$
Let $v_n\in M^n_{y_1^n,\ldots,y_{k_n}^n}$. 
By the definition of $\mu_{k_n}^n$,
\beq
\label{22.1}
I_n(v_n)\le \mu_{k_n}^n=I_n(u_n).
\eeq
By Taylor's formula,
\beq
\label{22.2}
\begin{array}{rcl}
\vspace{2mm}
I_{n} (v_n) - I_{n}(u_{n})
&=  & I'_{n}(u_{n})[v_n - u_{n}]
\\ \vspace{2mm}
& & + \frac{1}{2}
\int_{\R^N} (|\nabla(v_n - u_{n})|^2+ a(x)(v_n - u_{n})^2 )dx
\\
& & +\frac{1}{2} \int_{\R^N} (q  b_n(x) |\tau_n|^{q-1}- p  |\omega_n|^{p-1})(v_n -u_{n})^2 dx
\end{array}
\eeq
with $\tau_n (x) = u_{n}(x) + \tilde{\tau}_n (x)
(v_n -u_{n})(x)$ and $\omega_n (x) = u_{n}(x) + \tilde{\omega}_n (x)
(v_n -u_{n})(x)$, for suitable $\tilde{\tau}_n(x),\, \tilde{\omega}_n (x)\in [0,1]$.

Our aim is to obtain a contradiction with (\ref{22.1}) proving that a
careful estimate of the terms in (\ref{22.2}) implies that, for large $n$,
\beq
\label{38.1}
I_{n} (v_n)- I_{n}(u_{n}) > 0 .
\eeq
First step is proving that
\beq
\label{38.2}
\lim_{n\to  + \infty}|v_n- u_{n}|_\infty  = 0.
\eeq
By Proposition \ref{w} applied to $v_n$ and $u_{n}$,
considering that  $|y_1^n-x_1^n  |\lo 0$, we obtain
\beq
\label{38.3}
\lim_{n\to  + \infty} |v_n -u_{n}|_{L^\infty(\cup_{i=1}^{k_n}
B_{\bar{\rho}}(y_{i}^n))} = 0 \qquad \forall
\bar{\rho} > 0.
\eeq

On the other hand, on $\{0<u_n<\d\}\cap\{0<v_n<\d\}$ by Proposition \ref{P3.3}
$$
-\Delta u_n+[a_\infty-\alpha(x)]u_n+b_n(x)u_n^q=u_n^p,\qquad
-\Delta v_n+[a_\infty-\alpha(x)]v_n+b_n(x)v_n^q=v_n^p,
$$
hence, in view of (\ref{delta}), when $0< u_n< v_n<\d$ we deduce
\begin{eqnarray}
\nonumber
\Delta(v_n - u_{n}) & = &  (v_n - u_{n}) \left[(a_\infty-\alpha(x)) +b_n(x)
\frac{v_n^q -u_{n}^{q}}{v_n - u_{n}}
-\frac{v_n^p -u_{n}^{p}}{v_n - u_{n}} \right]  \\
\label{38.4}
& \geq & (a_0-p\d^{p-1}) (v_n - u_{n})\\
\nonumber
&\ge & 0.
\end{eqnarray}
Analogously, when $ 0<v_n < u_{n}  < \d $, we get
$$
\Delta(v_n -u_{n}) \leq 0.
$$
Hence, for $\bar{\rho}>R$,  we can say that  the maximum of
$|v_n - u_{n}|$   in
$\R^N \setminus \cup_{i=1}^{k_n}
B_{\bar{\rho}}(y_{i}^n)$ is attained on the boundary $\cup_{i=1}^{k_n}
\partial B_{\bar{\rho}}(y_{i}^n)$.
Then by (\ref{38.3}) we deduce
\beq
\label{39.2}
\lim_{n\to + \infty} |v_n -u_{n}|_{L^\infty(\R^N\setminus
  \cup_{i=1}^{k_n}  B_{\bar{\rho}}(y_{i}^n))} = 0,
\eeq
which together (\ref{38.3}) yields  (\ref{38.2}).

Set $ s_n = |v_n - u_{n}|_\infty$. 
Since $\beta_{x_1^n}(u_n)=\beta_{y_1^n}(v_n)=0$ and $x_1^n\neq y_1^n$,
$v_n \neq u_{n} $, and  $s_n>0$ so we can define
$$
\phi_n(x) = \frac{v_n(x) - u_{n}(x)}{s_n}.
$$
We set now  $\cI=\{1,\ldots,\bar k\}$, if $k_n\equiv\bar
k$, and  $\cI=\N$, if $\{k_n\}_n$ is unbonded, and
 $n(j)=\min\{n\in\N$ : $k_n\ge j \}$, for all $j  \in \cI$.

\vspace{2mm}

{\em  {\bf Claim} Up to a subsequence, for all  $j  \in \cI$, $\{\phi_n
  (x+ x_{j}^n)\}_{n\ge n(j )}$ converges in $H^1_{\loc}(\R^N)$ and
  uniformly on compact subsets of $\R^N$, to a solution $\phi$ of the equation
\beq
\label{*39}
- \Delta \phi + a_\infty\phi = p\, w^{p-1} \phi \qquad \mbox{ in }  \R^N.
\eeq
Moreover, the convergence is uniform with respect to $j\in\cI$.
}

\vspace{1mm}

We postpone the proof of this claim to the end of the argument.

Proposition \ref{w} implies
$\supp (v_n^\d)_1\subset B_{R}(  x_1^n)$, for large $n$,
hence we deduce
\beq
\label{39.1}
\int_{B_{R}(x_{i}^n)} [(v_n^\d )^2 - (u_{n}^\d)^2 ]
(x-x_{i}^n) dx  = \left\{
\begin{array}{cl}
0
&  \mbox{ if }i\in\{2,\ldots,k_n\}
\\
\sigma_n \bar \l \int_{B_{R}(x_{1}^n)} (v_n^\d)^2 dx
& \mbox{
  if }i=1
\end{array}
\right.
\eeq
because, when $i=1$:
\begin{eqnarray*}
\int_{B_{R}(x_{1}^n)} (v_n^\d )^2 (x- x_{1}^n
) dx & = &
\int_{B_{R}(x_{1}^n)} (v_n^\d)^2   (x-y_1^n)dx
+ \int_{B_{R}(x_{1}^n)}(v_n^\d)^2  ( y_1^n-
x_{1}^n ) dx\\
& = & \sigma _n  \bar\l \int_{B_{R}(x_{1}^n)} (v_n^\d)^2   dx.
\end{eqnarray*}
On the other hand, for all $i\in\{1,\ldots,k_n\}$
$$
\hspace{-4cm}\int_{B_{R}(x_{i}^n)} [(v_n)^\d )^2 - (u_{n}^\d)^2 ]
(x-x_{i}^n) dx=
$$
\beq
\label{40.1}
= 2  \int_{B_{R}(x_{i}^n)}  u_{n}^\d\, 
[v_n-u_{n}]
 (x-x_{i}^n)dx
 + \int_{B_{R}(x_{i}^n)} \cR _{v_n,u_{n}}(x)\,   (x-x_{i}^n)dx,
\eeq
where
$$
\cR_{v_n,u_{n}}  = (v_n^\d)^2
-(u_{n}^\d)^2 -
2 u_{n}^\d (v_n  -u_{n}).
$$
By a direct computation, using the convexity of
the real map
$t \mapsto g(t)=((t-\d)^+)^2$ and the fact that for all fixed $y \in
\R$ $\min_{t\in\R}[(t-y)^2 - \cR(t,y)]=0$ when $\cR(t,y) = g(t) - g(y)
- g'(y)(t-y)$, we deduce that 
\beq
\label{40*}
0 \leq \cR_{v_n,u_{n}}  \leq (v_n  -
u_{n})^2 \leq s_n^2.
\eeq
Comparing (\ref{39.1}) and (\ref{40.1}) and using (\ref{40*}),
we infer that when $i\neq 1$
\beq
\label{40.2}
\left|2 \int_{B_{R}(0)} (u_{n}(x+x_{i}^n))^\delta
 [v_n (x + x_{i}^n) - u_{n}(x+x_{i}^n)]\, x \,
 dx\right|\le c\,s_n^2,
\eeq
with $c$ independent of $i$, while when  $i=1$,
$$
2 \int_{B_{R}(0)} u_{n}^\delta (x+  x_{1}^n)
 [v_n (x +   x_1^n) - u_{n}(x+ x_1^n)]\,
x \, dx + O(s_n^2)
$$
\beq
\label{40.3}
 =\sigma_n  \bar\l \int_{B_{R}(0)} (v_n^\delta (x +  x_1^n))^2 dx.
\eeq
Therefore, since $s_n \neq 0$ and Lemma \ref{N2'} and (\ref{*39}) give
the existence of vectors $\tau_i \in \R^N$, $i \in \cI$, such that
\beq
\label{40.4}
\lim_{n\to  + \infty}\phi_n(x+x_{i}^n) = (\nabla w(x) \cdot
\tau_i)\qquad\forall i\in\cI,
\eeq
by using Proposition \ref{w}, equality (\ref{40.2}) and the
choice of $R$, we deduce
$$
0  = 2 \lim_ {n\to  + \infty}\int_{B_{R}(0)}
u_{n}^\d(x+x_i^n)  \phi_n (x + x_i^n)\, x \, dx
= 2 \int_{B_{R}(0)} w^\d (\nabla w \cdot \tau_i) \, x\, dx
$$
$$
= \int_{B_{R}(0)} (\nabla (w^\d)^2 \cdot \tau_i)\, x \, dx =
-  \tau_i \int_{B_{R}(0)} (w^\d)^2  dx\qquad\forall i\in \cI\setminus\{1\}
$$
which implies, for all  $i\in \cI\setminus\{1\}$, $\tau_i = 0$.
Analogously, when $i=1$, using (\ref{40.3}), we obtain
$$
- \tau_{1} \int_{B_{R}(0)} (w^\d)^2 dx=  \bar\l \left(\lim_{n\to  + \infty}\frac{\sigma_n }{s_n}
\right)\int_{B_{R}(0)}
(w^\d )^2  dx 
$$
which gives
\beq
\label{41.1}
\tau_{1} =
-   \left(\lim_{n\to  + \infty}\frac{\sigma_n }{s_n}\right)\, \bar\l.
\eeq
Now, we observe that $\tau_{1} \neq 0$.
In fact, otherwise, from (\ref{40.4}), we would derive that
$\lim_{n\to +\infty} |\phi_n|_{L^\infty(\cup_{i = 1}^{k_n}
  B_{\bar{\rho}}(y_{i}^n))} =
0$,
 for all $\bar{\rho} > 0$; moreover, by the argument
 used to obtain (\ref{39.2}),
  for $\bar{\rho}>R$  we could deduce the relation $\lim_{n\to  + \infty}$
  $|\phi_n|_{L^\infty(\R^N 
  \setminus\cup_{i = 1}^{k_n}
  B_{\bar{\rho}}(y_i^n))}
= 0$, and,
hence, $\lim_{n\to  + \infty} |\phi_n|_\infty = 0$,  contradicting
$|\phi_n|_\infty = 1$.
Thus, by (\ref{41.1}), $\sigma_n$ and $ s_n$ have the same order, namely
\beq
\label{42.2}
\lim_{n\to  + \infty} \frac{\sigma_n}{s_n} = \gamma \in \R^+ \setminus
\left\{0\right\}\quad\mbox{and} \quad \tau_{1} = -\gamma  \bar\l.
\eeq

Now, let $\tilde{\rho}\in (\frac{3}{4}R ,R)$ be fixed and
consider a cut-off decreasing function $\bar \chi \in C^\infty (\R^+,
[0,1])$ such that $\bar \chi(t) = 1$ if $t \leq \tilde{\rho}$, $\bar
\chi (t) = 0$ if $t \geq R$.
Thus, putting $\bar \chi_i^n (x) = \bar \chi (|x - y_{i}^n|)$,
with $i\in\cI$ and $n\ge n(i)$,  we have
$$
\bar \chi_i^n (x) [v_n (x) - u_{n}(x)] \in H^1_0
(B_{R} (y_{i}^n)),
$$
moreover, for large $n$, by Proposition \ref{w},
\beq
\label{42.1}
\left\{
\begin{array}{l}
\supp u_{n}^\delta  \subset
\cup_{i = 1}^{k_n}
  B_{\tilde{\rho}}(y_{i}^n)
\\
  \supp\ v_n^\d  \subset
\cup_{i = 1}^{k_n}
  B_{\tilde{\rho}}(y_{i}^n)
 \\
\bar\chi_i^n (x) = 1  \quad \forall x \in B_{\tilde{\rho}}(y_{i}^n) \ : \
u_{n}^\d(x) \neq 0.
\end{array}
\right.
\eeq

We are, now, in position of estimating the terms of the expansion
(\ref{22.2}).
Indeed, considering (\ref{42.1}), (\ref{eqs}),  (\ref{ML}), and
denoting by $\l^n_{y_{i}^n}$ the Lagrange multipliers related to $v_n$,
we can write, for large $n$:
\begin{eqnarray*}
 I'_{n}(u_{n})[v_n -u_{n} ]
& = &
I'_{n}(u_{n})
[(1 - \sum_{i = 1}^{k_n}\bar\chi_i^n) (v_n - u_{n})]+ I'_{n}(u_{n})
[\sum_{i = 1}^{k_n}\bar\chi_i^n(v_n - u_{n}) ]
\\
& = &
\sum_{i = 1}^{k_n}I'_{n}(u_{n})
[\bar\chi_i^n(v_n - u_{n}) ]
\\
& = &
\sum_{i = 1}^{k_n}
\int_{B_{R}(y_{i}^n)}u_{n}^\d\,
\bar\chi_i^n(v_n - u_{n})
 (\l^n_{y_{i}^n} \cdot (x- y_{i}^n)) dx
\\
& = &
\sum_{i = 1}^{k_n}
\int_{B_{R}(y_{i}^n)}u_{n}^\delta
(v_n - u_{n})
 (\l^n_{y_{i}^n} \cdot (x- y_{i}^n)) dx.
\end{eqnarray*}

Thus, considering (\ref{40.2}) and (\ref{42.2}), we obtain
$$
\lim_{n\to  + \infty}  \frac{1}{s_n |\l_1^n|}
 I'_{n}( u_{n})\left[v_n - u_{n} \right]
$$
\begin{eqnarray}
& = &\lim_{n\to  + \infty}\left[ \left( \frac{\l^n_1}{|\l^n_{1}|} \cdot
\int_{B_{R}(0)} u^\d_{n}(x+  x_1^n) \phi_n
(x +  x_1^n)  \, x \, dx \right) +
 \sum_{i =2}^{k_n} \frac{1}{s_n} \frac{\l^n_{y_{i}^n}}{|\l_1^n|}
O(s_n^2)
\right]
\nonumber \\
 & = &
-\gamma \int_{B_{R}(0)} w^\d(x) (\nabla
w(x)\cdot  \bar\l) ( \bar\l \cdot x)\, dx.
\label{44.1}
\end{eqnarray}

Moreover, taking into account (\ref{delta}), (\ref{54.2}) and (\ref{42.2}) we have
\begin{eqnarray}
&\liminf\limits_{n\to  + \infty}& \frac{1}{s_n |\l_1^n|}
  \int_{\R^N} \left\{ |\nabla (v_n - u_{n} )|^2 + a(x)(v_n
    - u_{n})^2 +\right.
\nonumber \\
&  & \left. +[q\,b_n(x)|\tau_n|^{q-1}-p \,|\omega_n|^{p-1}] (v_n - u_{n} )^2 \right\} dx
\nonumber \\
&\geq & \liminf_{n\to  + \infty}\frac{1}{s_n |\l_1^n|}
\left[\sum_{i = 1}^{k_n} \int_{B_{R}(y_{i}^n)} - p\,\omega_n^{p-1} (v_n - u_{n} )^2 dx \right.
 \nonumber \\
& & \hspace{1cm}\left. +\int_{\R^N\setminus
    \cup_{i=1}^{k_n}B_{R}(y_{i}^n)}
(a_0-p\d^{p-1}) (v_n - u_{n} )^2 dx\right]
\nonumber \\
&\geq & \lim_{n\to  + \infty} -C\,  \frac{ k_n s_n}{|\l_1^n|} = 0.
\label{44.2}
\end{eqnarray}
Finally, combining (\ref{22.2}), (\ref{44.1}) and (\ref{44.2}), we get
$$
 \liminf_{n\to  + \infty} \frac{I_{n} (v_n) - I_{n}(u_{n})}
{s_n |\l_1^n|}
 \geq
 -\gamma \int_{B_{R}(0)} w^\d (\nabla w \cdot  \bar\l)
( \bar\l \cdot x)dx > 0
$$
and as a consequence (\ref{38.1}), as desired.

\vspace{ 2mm}

 To complete the proof, let us now prove the claim.
Up to a diagonal argument, we can just prove the claim for a fixed $j\in\cI$.
Moreover, to simplify the notation, we use the same symbols to indicate
subsequences of a given sequence.

 Being
$u_{n} \in {M}^{n}_{  x_{1}^n,
  \ldots ,  x_{k_n}^n
}$
and
$v_n \in {M}^{n}_{
  y_{1}^n,
  \ldots ,  y_{k_n}^n }$, they verify respectively, by Corollary \ref{C1},  the
Euler-Lagrange equations
\beq
\label{delta1}
-\Delta u_{n} (x)+a(x) u_{n} (x)+b_n(x)u_n^q= u_{n}^{p} (x)
+ \sum_{i=1}^{k_n} (u_{n}^\delta)_i (x)(\lambda_i^n
\cdot (x-  x_{i}^n))
\eeq
\beq
\label{delta2}
-\Delta v_n (x)+a(x)v_n (x)+b_n(x)v_n^q=v_n^{p} (x)
 + \sum_{i=1}^{k_n}(v_n^\delta)_i(x) (\lambda_{y_i^n}^n
\cdot (x-  y_i^n)).
\eeq
Hence, taking into account $y_i^n=x_i^n$ for $i\in\{2,\ldots,k_n\}$,
we get 
$$
-\Delta(u_{n}-v_n)+ a(u_{n}-v_n)+b_n(x)(u_{n}^{q}-v_{n}^{q}) =(u_{n}^{p}-v_{n}^{p})
$$
$$
+(v_n^\delta)_1( \lambda_{y_1^n}^n\cdot(y_1^n-  x_{1}^n)) +(u_{n}^\d)_1((\lambda_1^n-
 \lambda^n_{y_1^n})\cdot(x-x_1^n))
+
[(u_{n}^\d)_1-(v_n^\d)_1]
(\lambda^n_{y_1^n}\cdot (x- x_1^n))
$$
\beq
\label{46.1}
+\sum_{i=2}^{k_n}(u_{n}^\d)_i((\lambda_i^n-
 \lambda^n_{y_i^n})\cdot(x-
x_{i}^n))+\sum_{i=2}^{k_n}
[(u_{n}^\d)_i-(v_n^\d)_i](\lambda_{y_i^n}^n\cdot (x- x_i^n)).
\eeq

Let us fix $j \in \cI$, for $n\ge n(j)$ set
$
\hat s_n=
\max\{s_n,|\lambda_i^n- \lambda^n_{y_i^n}|\}
$
and
$$
\hat\phi_n(x)= \frac{u_{n}(x +x_{j}^n )- v_n(x + x_{j}^n)}
{\hat{s}_n }.
$$
Remark that
\beq
\label{a_C}
|\hat{\phi}_n|_\infty \leq 1.
\eeq

Dividing by $\hat s_n$, we deduce from  (\ref{46.1})
$$
-\Delta\hat\phi_n(x)+a(x+x_{j}^n)\hat\phi_n(x)=
[\hat\omega_n (x+x_{j}^n )-b_n(x+x_{j}^n )\hat\tau_n(x+x_{j}^n )]\hat\phi_n(x)
$$
$$
+(v_n^\d)_1(x+x_{j}^n)
\left(\hspace{-1mm} \lambda^n_{y_1^n}\cdot {y_1^n-  x_{1}^n\over
    \hat s_n}\hspace{-1mm}\right)+(u_{n}^\d)_1(x+x_{j}^n)
\left(\hspace{-1mm}
{\lambda_1^n-  \lambda^n_{y_1^n}\over \hat s_n}
\cdot(x+x_{j}^n-  x_1^n)
\hspace{-1mm}\right)
$$
$$
+{(u_{n}^\delta)_1(x+x_{j}^n)
  -(v_n^\d)_1(x+x_{j}^n )\over \hat s_n}
\,
(\lambda^n_{y_1^n}\cdot (x+x_{j}^n- x_1^n))
$$
$$
+\sum_{i=2}^{k_n}(u_n^\d)_i(x+x_j^n)\left(
{\lambda_i^n-
 \lambda^n_{y_i^n}\over \hat s_n}\cdot(x+x_j^n-x_i^n)\right)
$$
\beq
\label{f}
+\sum_{i=2}^{k_n}
{(u_{n}^\d)_i(x+x_j^n)
  -(v_n^\d)_i(x+x_j^n)
\over   \hat s_n}
\,
(\lambda_{y_i^n}^n\cdot (x+x_j^n- x_i^n))
\eeq
where

$\hat\omega_n(x)={u_{n}^p(x)-v_n^p(x)\over
u_{n}(x)-v_n(x)}=p\,\check\omega_n^{p-1}(x)$ and $\hat\nu_n(x)={u_{n}^q(x)-v_n^q(x)\over
u_{n}(x)-v_n(x)}=q\,\check\nu_n^{q-1}(x)$,
$\check\omega_n$ and $\check\nu_n$ being  functions taking their
values between the values 
of $u_{n}$ and $v_n$. 
We also remark that the relation
\beq
\label{b_C}
\frac{|\l_j^n -{\l}^n_{y_j^n}
  |}{\hat{s}_n} \leq 1 \quad\mbox{ if }j \in \cI
\eeq
clearly hold true, while the relations
\beq
\label{b2}
\l_j^n\lo 0, \quad
 \l^n_{y_j^n}\lo 0\quad \mbox{ if }j \in \cI
\eeq
hold uniformly with respect to the choice of $j$ and can be obtained
arguing as for proving (\ref{5.17}), by
using (\ref{delta1}) and (\ref{delta2}) respectively.
Moreover, since for large $n$ the equality
$$
y_1^n-  x_1^n =
\frac{1}{|(v_n^\d)_1|_2^2} \int_{B_{R}(x_1^n)}\hspace{-2mm} x\ [(v_n^\d)_1 (x)]^2 dx -
\frac{1}{|(u_{n}^\d)_1|_2^2} \int_{B_{R}(x_1^n)} \hspace{-2mm}x \
[(u_{n}^\d)_1 (x)]^2 dx
$$
holds, considering that
$|(v_n^\d)_1 (x)-(u_{n}^\d)_1 (x)| \leq s_n$,
 $(v_n^\d)_1(\cdot+y_1^n) \lo w^\d$ and
$(u_{n}^\d)_1 $ $(\cdot + x_1^n)\lo w^\d$,
 by a direct computation, we deduce
\beq
\label{b'}
|y_1^n - x_1^n| \leq  c   s_n,\qquad\mbox{ for large }n.
\eeq

Now, let us observe that by (\ref{f})
$$
 \hspace{-1cm}- \Delta \hat{\phi}_n + a(x +x_j^n) \hat{\phi}_n(x)
=[\hat\omega_n (x+x_{j}^n )-b_n(x+x_j^n)\hat\nu_n(x+x_{j}^n )]
\hat{\phi}_n (x ),
$$
$$ \hspace{7cm} \forall
x\in\R^N\setminus\cup_{i=1}^{k_n}B_{R+\sigma_n|\bar\lambda|}(x_i^n-x_j^n),
$$
where $|\hat\omega_n | \leq p\, \d^{p-1}$ and  $|\hat\nu_n | \leq q
\, \d^{q-1}$. 
Therefore $|\hat{\phi}_n |$ satisfies
$$
-\Delta|\hat{\phi}_n (x)| + (a_0 - p\, \d^{p-1})|\hat{\phi}_n|\leq 0\quad\mbox{ in }
\R^N\setminus\cup_{i=1}^{k_n}
B_{R+\sigma_n|\bar\lambda|}(x_i^n-x_j^n).
$$
Then, by Lemma \ref{N*}, choosing $\hat
r>R+\sigma_n|\bar\lambda|$ and taking into account Proposition
\ref{P4.1}, we obtain $|\hat\phi_n(x)| < c e^{-\xi |x|}$, $x\in
B_{2\hat r}(0)\setminus  B_{R + \sigma_n |\bar\l|} (0)$, for $\xi \in
(0, \sqrt{a_0  - p \d^{p-1}})$, $c$ independent of $\hat r$ and large $n$.
Therefore, taking into account (\ref{a_C}), we deduce the existence of
$C>0$ such that $\forall \hat r>0$ and for large $n$
\beq
\label{g}
\int_{B_{2 \hat{r}} (0)}|\hat{\phi}_n| dx \leq C.
\eeq
Thus,
using (\ref{a_C}), (\ref{g}), (\ref{b_C}), (\ref{b2}), (\ref{b'}) and
observing that, by Proposition \ref{w},
$u_{n}(x+x_j^n)\lo w (x)$,
$v_n(x+x_j^n)\lo w (x)$,
$|\check{\omega}_n (x+x_j^n)|\lo w(x)$ and $|\check{\nu}_n (x+x_j^n)|\lo w(x)$,
 from (\ref{f}) we deduce for large $n$
\beq
\label{h}
\int_{B_{2 \hat{r}} (0)} | \Delta \hat{\phi}_n||\hat{\phi }_n| dx \leq
C_1,
\eeq
where $C_1$ is a real positive constant depending neither on $n$
nor on $\hat{r}$.
Fix now a cut-off function $\tilde{\xi} \in \cC^\infty_0(B_2(0),[0,1])$ such that
$\tilde{\xi} = 1$ on $B_1(0)$. 
Setting
$\tilde{\xi}_{\hat{r}}(x) = \tilde{\xi} (\frac{x}{\hat{r}})$,
remark that $\int_{B_{2\hat{r}}(0)}|\Delta \tilde{\xi}_{\hat{r}}|
= C_2 \hat{r}^{N-2}$.
Then, using (\ref{h}), we get the relation
\begin{eqnarray*}
\int_{B_{\hat{r}}(0)}|\nabla \hat{\phi}_n |^2
& \leq &
\int_{\R^N} |\nabla \hat{\phi}_n |^2 \tilde{\xi}_{\hat{r}}dx
 =
- \int_{\R^N}
\hat{\phi}_n \nabla ( \tilde{\xi}_{\hat{r}}\nabla \hat{\phi}_n )
dx
\nonumber \\
& = &
- \int_{\R^N} \hat{\phi}_n \tilde{\xi}_{\hat{r}} \Delta\hat{\phi}_n dx
- \frac{1}{2}\int_{\R^N}
\nabla \hat{\phi}_n^2 \nabla \tilde{\xi}_{\hat{r}}
dx
\nonumber \\
& \leq &
\int_{B_{2 \hat{r}} (0)} | \Delta \hat{\phi}_n| |\hat{\phi}_n|dx
+ \frac{1}{2}
\int_{B_{2 \hat{r}} (0)} \hat{\phi}_n^2 \Delta \tilde{\xi}_{\hat{r}}
dx
\nonumber\\
& \leq & C_1+{C_3\over 2}\, e^{- 2 b \hat{r}} \hat{r}^{N-2}
\end{eqnarray*}
which implies that $(\hat{\phi}_n)_n$ is bounded in $H^1_{\loc}(\R^N)$.

Hence, we can assume that a function $\phi$ exists such that,
up to a subsequence, $\hat{\phi}_n\rightharpoonup \phi $ in
$H^1_{loc}(\R^N)$ and, in view of (\ref{a_C}),
$|\hat{\phi}_n - \phi|_{L^q(K)} \to 0$ for all $q <  + \infty$ and
for all compact sets $K\subset\R^N$.
Furthermore, we can pass to the limit in (\ref{f}) and  obtain
\beq
\label{l}
-\Delta \phi +a_\infty \phi - p w^{p-1} \phi = w^\d (\lambda' \cdot x),
\eeq
where  $\lambda' = \lim_{n\to  + \infty }\frac{\lambda_j^n -{\lambda}^n_{y_j^n}}{\hat{s}_n}$, for
  $j \in \cI$.

Now, to complete the argument, we only need to show that $\l' = 0$.
Indeed, in this case, for large $n$, $\hat{s}_n = s_n$,  $\hat{\phi}_n
= \phi_n$ and what above proved for $\hat{\phi}_n$ is just what
asserted in the claim.
From (\ref{l}), using Lemma \ref{N2'} and Fredholm alternative
theorem, we deduce that $w^\d (\l'\cdot x)$ must be orthogonal to
$\frac{\partial w}{\partial \nu}$, for all $\nu\in\R^N\setminus\{0\}$,  so,
choosing $\nu = \l'$, we obtain
\begin{eqnarray*}
0& =& \int_{\R^N} \frac{\partial w}{\partial \l'}\,  w^\d(\l'\cdot
x)dx = \frac{1}{2}\int_{\R^N} \frac{\partial (w^\d)^2}{\partial\l'}(\l'\cdot x)dx
 \\
& = & -\frac{1}{2}\int_{\R^N} (w^\d)^2 \frac{\partial
}{\partial\l'} (\l'\cdot x)dx = -\frac{1}{2} |\l'|^2\int{\R^N}  (w^\d)^2 dx,
\end{eqnarray*}
that gives $\l' = 0,$ as desired.

\qed

\vspace{1cm}


\appendix
\makeatletter 
\newcommand{\section@cntformat}{Appendix \thesection:\ }
\makeatother

\section{Appendix: On the existence of Ground States} 
\setcounter{equation}{0}

As observed in the introduction, a natural question to wonder is
either in our assumptions problem $(P)$ admits a ground state solution
or all the solutions of $(P)$ are bound states.
Here we show, by means of some examples, that both the situations can
be true. 
We refer the reader also to \cite{CeP} for a deeper study of the same
question with respect to another equation with competing coefficients.
Clearly in what follows we shall be sketchy on some well known facts.

\vspace{5mm}

\noindent {\em {A case in which a ground state solution exists}} 

\vspace{2mm}

\noindent We consider $(P)$ with coefficients $a(x)$ and $b(x)$ chosen as follows
\beq
\label{E1}
a(x)=2-e^{-|x|},\qquad b(x)={C\over 1+|x|}.
\eeq
Clearly $C$ can be chosen so small that all the assumptions of Theorem
\ref{T} are fulfilled and the theorem applies.

Now, setting
$$
m_a:=\inf\left\{\int_{\R^N}(|\D u|^2+a(x)u^2)\, dx\ :\ u\in H^1(\R^N),\
\int_{\R^N}|u|^{p+1}=1\right\},
$$
and testing $\int_{\R^N}(|\D u|^2+a(x)u^2)dx$ by ${w(x)\over
  |w|_{L^{p+1}}}$ ($w$ being the function defined in Lemma \ref{N2})
  one easily obtains
$$
m_a<\min \left\{\int_{\R^N}(|\D u|^2+2 u^2)\, dx\ :\ u\in H^1(\R^N),\
\int_{\R^N}|u|^{p+1}=1\right\}
$$
and a standard application of the Concentration-Compactness principle
allows to conclude that $m_a$ is achieved.
Then, it is clear that a function realizing 
$$
\inf\left\{\int_{\R^N}(|\D u|^2+a(x)u^2)\, dx+\int_{\R^N}b(x)\,
  |u|^{q+1}dx\ :\ u\in H^1(\R^N),\ \int_{\R^N}|u|^{p+1}=1\right\}
$$
must exist too, by continuity, when $C$ is suitably small.
Therefore, for small $C$, $(P)$ has a ground state solution and
infinitely many positive bound state solutions.

\vspace{5mm}

\noindent {\em{A case in which there is no ground state}} 

\vspace{2mm}

\noindent We consider, for all $n\in\N$, $(P)$ with coefficients
\beq
\label{E2}
a_n(x)=1-{1\over 2} \chi_{B_{1/n}(0)}(x),\qquad b(x)={C\over 1+|x|}.
\eeq

For all $n$, if $C$ is suitably small, all the assumptions of Theorem
\ref{T} are satisfied so its conclusions hold true.
Furthermore the arguments used to prove Theorem \ref{T} show that a
constant $\bar C>0$, not depending on $n$, exists so that when $a_n$
and $b$ are as in (\ref{E2}) and $C\in(0,\bar C)$, $(P)$ has infinitely many positive solutions whatever $n\in\N$
is.

In what follows we assume $C\in(0,\bar C)$ fixed and we intend to show
that $(P)$ cannot have ground state solutions for large values of $n$.

To carry out our argument we first need to introduce some notation.
We set
$$
I_n(u):={1\over 2}\int_{\R^N}[|\D u|^2+a_n(x)u^2]dx+{1\over
  q+1}\int_{\R^N}b(x)|u|^{q+1}dx-{1\over
  p+1}\int_{\R^N}|u|^{p+1},
$$
$$
\cN_n:=\{u\in H^1(\R^N)\ :\ I'_n(u)[u]=0\},
\qquad
\cN_\infty:=\{u\in H^1(\R^N)\ :\ I'_\infty(u)[u]=0\},
$$
where $I_\infty$ is the functional related to the limit problem
and defined in Section 2.
Standard arguments, similar but simpler than those displayed in Lemma
\ref{lMax}, show that
\begin{itemize}
\item
$\forall u\in H^1(\R^N)\setminus\{0\}\ \exists$ unique
$t_u\in(0,+\infty)$ s.t. $t_uu\in\cN_\infty, \
I_\infty(t_uu)\hspace{-1mm}=\hspace{-2mm}\max\limits_{t\in(0,+\infty)}I_\infty(tu)$
\item
$\forall u\in H^1(\R^N)\setminus\{0\}\ \exists$ unique
$\tau_u^n\in(0,+\infty)$ s.t. $\tau_u^nu\in\cN_n, \
I_n(\tau_u^nu)\hspace{-1mm}=\hspace{-2mm}\max\limits_{\tau\in(0,+\infty)}I_n(\tau u)$
\item
$u\in\cN_n\ \Rightarrow\ \tau^n_u=1,\ 
u\in\cN_\infty\ \Rightarrow\ t_u=1. 
$
\end{itemize}

Moreover, for all $n\in\N$, for all $u\in\cN_n$ we have
\begin{eqnarray}
\label{16.50} 0&=& \int_{\R^N}(|\D u|^2dx+a_nu^2)dx+\int_{\R^N}b(x)\,
|u|^{q+1}dx-|u|^{p+1}_{p+1}\\
&\ge & {1\over 2}\|u\|^2-|u|^{p+1}_{p+1}\ge {1\over 2}\|u\|^2-c\|u\|^{p+1},\nonumber
\end{eqnarray}
from which
\beq
\label{19.24}
\|u\|\ge \tilde c,\qquad 
|u|_{p+1}\ge\tilde c \qquad \forall u\in\cN_n,
\eeq
with $\tilde c$ constant independent of $n$.

Furthermore we remark that the inequality
$$
m_n:=\inf_{\cN_n}I_n(u)\le\inf_{\cN_\infty}(u)=:m_\infty
$$
holds true $\forall n\in\N$, as one can verify testing $I_n$ by a
sequence $\tau^n_{w_n}w_n(x)$, with $w_n(x)=w(x-y_n)$,
$|y_n|\to+\infty$.

Now, let us argue by contradiction and let us assume that for all
$n\in\N$ a ground state solution $u_n$ of $(P)$, with coefficients
$a_n$ and $b$, exists.
We can suppose $u_n\ge 0$, because otherwise we can replace it by
$|u_n|$, then the maximum principles gives $u_n>0$.

We claim that
\beq
\label{o}
u_n(x)=w(x-y_n)+\phi_n(x)
\eeq
where $|y_n|\to+\infty$ and $\phi_n\in H^1(\R^N)$, $\phi_n\to 0$
strongly in $H^1(\R^N)$.

Indeed from
\begin{eqnarray*}
m_\infty \ge I_n(u_n)& = & \left({1\over 2}-{1\over
    p+1}\right)\|u_n\|^2-{1\over 2}\left({1\over 2}-{1\over
  p+1}\right)\int_{B_{1/n}(0)} u_n^2dx\\
& & +\left({1\over q+1}-{1\over
    p+1}\right)\int_{\R^N}b(x)\, u_n^{q+1}dx\\
& > & {1\over 2}\left({1\over 2}-{1\over
  p+1}\right)\|u_n\|^2,
\end{eqnarray*}
we deduce that $\{\|u_n\|\}_n$ is bounded, from which, setting
$t_n:=t^n_{u_n}$ and taking into account that
$t_n=\left({\|u_n\|^2\over |u_n|_{p+1}^{p+1}}\right)^{1\over p-1}$, we
also obtain $\{t_n\}_n$ bounded.
Thus
\begin{eqnarray*}
m_\infty &\le & I_\infty(t_nu_n)\le{1\over 2}\|t_nu_n\|^2-{1\over
  p+1}|t_nu_n|^{p+1}_{p+1}+{1\over q+1}\int_{\R^N}b(x)\, (t_nu_n)^{q+1}dx\\
& = & I_n(t_nu_n)+o(1)\le I_n(u_n)+o(1)\\
& \le & m_\infty+o(1),
\end{eqnarray*}
that implies 
\beq
\label{10.37}
\lim_{n\to+\infty}I_\infty(t_nu_n)=\lim_{n\to+\infty}I_n(t_nu_n)=\lim_{n\to+\infty}I_n(u_n)=
m_\infty.
\eeq 
So, being $t_nu_n\in \cN_\infty$, Lemma \ref{N2} implies 
\beq
\label{17.35}
t_nu_n=w(\cdot-y_n)+\phi_n,
\eeq
with $y_n\in\R^N$ and $\phi_n\to 0$ strongly in $H^1(\R^N)$, and, by 
Schauder estimates (see, f.i., \cite{serrin}), in $\cC^2_{\loc}(\R^N)$
too.
Now, in view of (\ref{19.24}), we have $\lim_{n\to +\infty}t_n=\bar
t>0$ up to a subsequence, and, more, $\bar t=1$, because otherwise the
impossible inequality
$$
m_\infty =\lim_{n\to +\infty}I_\infty(t_nu_n)=I_\infty(\bar tw)<I_\infty(w)=m_\infty
$$ 
had to be true.
Finally $|y_n|\to +\infty$ is obtained because $\lim_{n\to
  +\infty}y_n=\bar y\in\R^N$ would imply
\begin{eqnarray*}
\lim\limits_{n\to +\infty}I_n(t_nu_n) & = &  \lim\limits_{n\to+\infty}
\left[I_\infty(t_nu_n)+{1\over 4}\int_{B_{1/n}(0)}u_n^2dx+
{1\over q+1}\int_{\R^N}b(x)\, u^{q+1}_ndx\right]\\
&=& m_\infty+{1\over q+1}\int_{\R^n}b(x)\, w^{q+1}(x-\bar
y)\, dx \ > \ m_\infty,\\
\end{eqnarray*}
contradicting $I_n(t_nu_n)\le I_n(u_n)\le m_\infty$.

Now, relation (\ref{o}) allows us to assert
\beq
\label{17.33}
\int_{\R^N}b(x)\, (t_nu_n)^{q+1}dx\ge c\,{1\over |y_n|}\qquad\forall n\in\N.
\eeq
On the other hand, being $I_\infty(t_nu_n)<I_n(t_nu_n)$  impossible
because it would imply
$$
m_\infty\le I_\infty(t_nu_n)< I_n(t_nu_n)\le I_n(u_n)\le m_\infty,
$$
we conclude that the opposite relation
\beq
\label{17.23}
I_n(t_nu_n)\le I_\infty(t_nu_n)
\eeq
must be true.
Therefore
\beq
\label{essa}
{1\over q+1}\int_{\R^N}b(x)(t_nu_n)^{q+1} dx
-{1\over 4}\int_{B_{1/n}(0)}( t_nu_n)^2\le 0
\eeq
follows and, in view of (\ref{17.33}), we can get the desired
contradiction if we show that 
\beq
\label{17.37}
u_n(x)\le c\, e^{-\sigma  |y_n|} \ \mbox{  on } B_{1/n}(0)
\eeq
holds for some $\sigma>0$, with $c$ independent of $n$.
To this end put $v_n(x)=u_n(x+y_n)$ and remark that to get
(\ref{17.37}) it is enough to show 
\beq
\label{xo}
v_n(x)\le c\, e^{-\sigma  \dist (x,B_R(0))}.
\eeq
To show (\ref{xo}) we start observing that for all $n\in\N$, $v_n$ solves
$-\Delta u+a_n(x+y_n) u +b(x+y_n)u^q -u^p = 0$ in $\R^N\setminus
B_{R}(0)$.
Then, choosing $\d$ as in (\ref{delta}) and considering that  $\inf a_n=1/2$ for all
$n\in\N$, to obtain (\ref{xo}) we need to show
\beq
\label{9.17}
v_n(x)<\d\qquad  \mbox{ on }\R^N\setminus B_R(0).
\eeq
Indeed, if (\ref{9.17}) is true, for all $n\in\N$ $v_n$ turns out to
be solution of 
$$
\left\{\begin{array}{lc}
-\Delta u+\left({1\over 2}-\d^{p-1}\right)u\le 0 &\mbox{ in
}\R^N\setminus B_R(0)\\
0<u_\d\le\d &\mbox{ in }\R^N\setminus B_R(0),
\end{array}\right.
$$
and (\ref{xo}) comes as consequence of Lemma \ref{N*}.

Therefore, what is left to complete the argument is to prove
(\ref{9.17}).
Again we argue by contradiction, so we assume the existence of a
sequence $\{z_n\}_n$, $z_n\in\R^N$, such that, up to a subsequence,
\beq
\label{18.40}
\{x\in B_{R/2}(z_n)\ :\ |x|>R,\  v_n(x)\ge \d\}\neq\emptyset\qquad \forall
n\in \N.
\eeq
Since $v_n\to w$ in $\cC^2_{\loc}(\R^N)$ and $w<\d$ on $\R^N\setminus
B_{R/2}(0)$, we deduce $|z_n|\to +\infty$ and, in turn, $v_n\to 0$ in
$H^1(B_{R/2}(z_n))$.
Then, being $v_n$ solution of $-\Delta u+a_n(x+y_n) u +b(x+y_n)u^q
-u^p = 0$ in $B_{R/2}(z_n)$, regularity arguments imply $v_n\to 0$
uniformly in $B_{R/2}(z_n)$ and, hence, give a contradiction with (\ref{18.40}).

\vspace{1cm}

{\small {\bf Acknowledgement}. The authors have been supported by the ``Gruppo
Nazionale per l'Analisi Matematica, la Probabilit\`a e le loro
Applicazioni (GNAMPA)'' of the {\em Istituto Nazionale di Alta Matematica
(INdAM)} - Project: Sistemi differenziali ellittici nonlineari
derivanti dallo studio di fenomeni elettromagnetici.

The second author has been also supported by the 2015 - Italian PRIN
Project number 2015KB9WPT.}

\end{document}